\documentclass[11pt,a4paper]{article}

\usepackage[english]{babel}
\usepackage{verbatim}
\usepackage{bm}

\usepackage{amsmath}
\usepackage{amssymb}
\usepackage{amsfonts}
\usepackage{amsbsy}
\usepackage{xcolor}
\usepackage{graphicx}
\usepackage{graphics}
\usepackage{subfig}
\usepackage{epsfig}
\usepackage{cite}
\usepackage{array}
\usepackage{booktabs} 
\usepackage{caption}  
\usepackage{tabularx}

\usepackage{hyphenat} 

\usepackage{amsthm} 

\newtheorem*{proof*}{Proof}
\newtheorem*{remark*}{Remark}

\numberwithin{equation}{section}
\numberwithin{theorem}{section}
\numberwithin{lemma}{section}
\numberwithin{defin}{section}
\numberwithin{cor}{section}
\numberwithin{prop}{section}

\begin{document}
\title{   Multi-rate Runge-Kutta methods:\\ stability analysis and applications}
\author{Bernhard Bachmann$^{(1)}$, Luca Bonaventura$^{(2)}$, Francesco Casella$^{(3)}$ \\
Soledad Fern\'andez-Garc{\'\i}a$^{(4)}$, Macarena G\'omez-M\'armol$^{(5)}$\\
{    Philip Hannebohm $^{(1)}$}
 }
\maketitle
\begin{center}
{\small
$^{(1)}$ Faculty of Engineering and Mathematics\\
 University of Applied Sciences and Arts, Bielefeld, Germany\\
 {\tt bernhard.bachmann@hsbi.de}, {\tt philip.hannebohm@hsbi.de}
 }
\end{center}
\begin{center}
{\small
$^{(2)}$  Dipartimento di Matematica\\
 Politecnico di Milano \\
{\tt luca.bonaventura@polimi.it}
}
\end{center}

\begin{center}
{\small
$^{(3)}$  Dipartimento di Elettronica, Informazione e Bioingegneria\\
 Politecnico di Milano \\
{\tt casella@elet.polimi.it}
}
\end{center}
\begin{center}
{\small
$^{(4)}$  
Departamento
de Ecuaciones Diferenciales y An\'alisis Num\'erico \& IMUS, \\
Universidad de Sevilla\\
{\tt soledad@us.es}
}
\end{center}
\begin{center}
{\small
$^{(5)}$
Departamento
de Ecuaciones Diferenciales y An\'alisis Num\'erico, \\
Universidad de Sevilla\\
{\tt macarena@us.es}
}

\end{center}

\date{}


\vspace*{1.0cm}


\vspace*{1.0cm}

\pagebreak

\abstract{We present an approach for the efficient implementation of self-adjusting multi-rate Runge-Kutta methods and we {   introduce a novel stability analysis, that covers the multi-rate extensions of all standard Runge-Kutta methods and allows to assess the impact of different interpolation methods for the latent variables and of the use of an arbitrary number of sub-steps for the active variables. The stability analysis applies successfully to the model problem typically used in the literature for multi-rate methods. Furthermore,} we also propose a physically motivated model problem that can be used to assess stability {   to problems with purely imaginary eigenvalues and in situations closer to those arising in applications.}
Finally, we present {   an efficient implementation
of multi-rate Runge-Kutta methods in the framework of the \textit{OpenModelica} open-source modelling and simulation software. }
Results of several numerical experiments,  performed with this implementation of the proposed methods, demonstrate the efficiency gains deriving from the use of the proposed multi-rate approach for physical modelling problems with multiple time scales.}

\pagebreak

\section{Introduction}
 \label{sec:intro} \indent
 We consider numerical methods for systems of Ordinary Differential Equations (ODEs) that
can be partitioned into components $y_s, y_f$ with slow and fast dynamics, respectively, as 
   \begin{equation}
   \label{eq:partitioned_sys}
   {  y}^{\prime} =\left[\begin{array}{ll}
   {  y}^{\prime}_{s} \\
   {  y}^{\prime}_{f}
\end{array} 
\right ]=
\left[\begin{array}{ll}
  f_s(y_s, y_f, t) \\
  f_f(y_s, y_f, t)
\end{array} 
\right ] = f(y,t)
\end{equation}
for $t\in[0,T]$.
 For such systems, several methods have been proposed {  
 over the last 60 years},
 in which different time steps are employed for the slow and fast components,   {    see e.g.  \cite{andrus:1979,gear:1984,rice:1960}.}
 These  methods
are generically known as multi-rate methods, as opposed to standard single-rate methods for ODEs. 
Variables $y_s, y_f $ are also known in the literature as latent and active variables, respectively.
 More recently, a number of methods have been proposed that do not require \textit{a priori} partitioning as in   \eqref{eq:partitioned_sys}, 
 but rather allow to identify fast and slow variables at runtime by their compliance, or lack of compliance, with a given error indicator, so that
the system has   form \eqref{eq:partitioned_sys} 
with a  right-hand side partition
that can be different for each time discretization interval. A review of these methods is presented in \cite{bonaventura:2020a}, where a variant of this more general, self-adjusting multi-rate approach was introduced.
The proposal in \cite{bonaventura:2020a} was tailored on
the specific implicit TR-BDF2 method  \cite{bank:1985, hosea:1996}, along  with a  general stability analysis for one-step multi-rate methods.

{   The present work has three main goals.
Firstly, we propose a
novel approach that combines
the self-adjusting multi-rate technique with standard time step adaptation methods.
}
By marking as fast variables only a small percentage of the variables associated with the
largest values of an error estimator, we obtain a more effective multi-rate implementation of a generic Runge-Kutta (RK) method. 
 Whenever the global time step is sufficient to guarantee a given error tolerance for the slow variables, but not for the fast ones, the multi-rate procedure is employed to achieve uniform accuracy at reduced computational cost. This approach is similar to that used successfully in spatially adaptive finite element techniques, see e.g. \cite{arndt:2022,
bangerth:2007, orlando:2022a, orlando:2022b},
where a fixed percentage of the simulation Degrees of Freedom (DOF) can be marked for refinement
if error indicator values exceed a given tolerance. 
{  
We then review the literature on the stability analysis of multi-rate methods and highlight
the lack of a general analysis
for multi-rate RK methods. We present a stability analysis for a general multi-rate RK method, that allows to study linear stability for a generic interpolation procedure and for
an arbitrary number of sub-steps for the active components.
The analysis covers the cases of Explicit RK methods (ERK), Diagonally Implicit RK methods (DIRK) and Singly Diagonally Implicit Runge-Kutta methods with Explicit first stage (ESDIRK). 
As discussed in the detailed review of the analyses of multi-rate methods available in the literature, a linear stability analysis in this very general framework does not appear to have been presented so far.
We apply the analysis results to the typical 2-DOF model problem used in the literature to assess the stability of multi-rate extensions of classical RK methods and the impact of different interpolation methods for the latent variables. 
We show that the multi-rate extension does not entail any loss of stability in this case. Indeed, for problems with fast variables only mildly coupled to the slow ones, the stability regions of the standard single rate RK methods are effectively increased in size by the multi-rate extension.
}

In order to apply this methodology {   also }to a case that is {  more }relevant
for applications, we {   then propose a more physically motivated  4-DOF model problem. Also in this case, we verify that, } in regimes of
mild coupling between fast and slow variables and in the presence of sufficient dissipation, multi-rate approaches generally maintain the stability properties of the corresponding single-rate methods,{   while a loss of unconditional stability of implicit methods arises in the case of stronger coupling and purely imaginary eigenvalues.}

{   Finally, we present the first implementation of multi-rate extensions of classical RK methods in an advanced open-source software package for industrial simulations. More specifically, we have implemented the proposed multi-rate approach in the framework of the OpenModelica software \cite{fritzson:2020}. We present the results of the application of this implementation to a number of numerical benchmark problems.}
 Due to their computational advantages, we will focus in particular on the application of multi-rate versions of higher order ESDIRK methods.
A very comprehensive review of these methods was presented in \cite{kennedy:2016, kennedy:2019},
where the properties and potential advantages of this class of methods are discussed in detail. 
The results demonstrate the efficiency gains deriving from the use of the proposed multi-rate approach for problems with multiple time scales.
{   Furthermore, the results of the numerical experiments highlight how different optimization strategies should be applied in the multi-rate and single-rate case, in order to take full advantage of the multi-rate DOF reduction.}

The outline of the paper is the following.
In Section \ref{sec:multirate}, the self-adjusting multi-rate approach is presented.
{   In Section \ref{sec:review}, a thorough review of the literature on multi-rate stability analysis is presented.}
 In Section \ref{sec:stability}, a linear stability analysis of the proposed methods is performed, which considers an arbitrary number of sub-steps for the active variables.
 In Section \ref{sec:model}, a specific model problem is considered, to which the previously introduced analysis template is applied. The resulting analysis highlights the impact on multi-rate stability of the
different choices for the interpolation operator of
the latent variables and of the number of sub-steps for the
active variables. {   In Section \ref{sec:openmod}, we briefly introduce the OpenModelica software package and describe the implementation of multi-rate methods in that context. }  In Section \ref{sec:tests}, numerical results are presented, which demonstrate the good performance of the proposed methods
even against optimal implementations of single-rate solvers.
  In Section \ref{sec:conclu}, we summarize our results and discuss possible
future developments.

\section{Effective implementation of multi-rate methods}
 \label{sec:multirate} 
 
 \indent
 
 Multi-rate methods have been defined in a number of previous papers. Here, we follow the general
outline of \cite{bonaventura:2020a}, but we introduce several simplifications and extensions which allow to achieve a more
efficient implementation.
 Consider the generic nonlinear system
$ y^{\prime} = f(y,t)$,  where $y: [0,T]\rightarrow \mathbb{R}^N $ is the solution of the continuous first order ODE problem
and $f:\mathbb{R}^N\times [0,T]\rightarrow \mathbb{R}^N $  the ODE right-hand side, which is assumed to satisfy the usual
regularity requirements to guarantee local existence and smoothness of the solutions.
Let then $t_n, n\in \mathbb{N} $ denote a set of discrete time levels such that
$t_0=0$, $t_{n+1}=t_n+h_n$ and $h_n$ denotes a generic time step.
At each time level $n,$ two approximations $u_{n+1}$ and $\hat u_{n+1}$ of the solution $y(t_{n+1})$ are computed,  with convergence orders $p$ and $\hat p$, respectively. We set $q = \min (p , \hat p). $ In the case of RK methods,  approximation
 $\hat u_{n+1}$  is usually obtained from an embedded method,
see e.g. \cite{hairer:2000}, to which we refer for all the general concepts on ODE methods used in the following.
An error estimator for the less accurate approximation is
then given by $ \|u_{n+1} - \hat u_{n+1}\|$.
 
 In order to comply with assigned error tolerances
  $\tau^r, $ $\tau^a$   for relative and absolute errors, respectively, standard single-rate implementations of RK methods introduce for each time step $i=1,\dots, N$ the quotient
 \begin{equation}
 \label{eq:stepSizeControl_01}
	  \eta^{n+1}_{i}= \frac{| u^{n+1}_{i} - \hat u^{n+1}_{i} | }{ \tau^r | u^{n+1}_{i} | + \tau^a} 
 \end{equation}
  and   require that the inequality
 $$\eta =  \max_{i=1,\dots,N} \eta_{i}^{n+1} \leq 1 $$ 
 holds. 
An optimal choice of the new time step value for an adaptive
single-rate implementation is then given by
$$h_{new} = h_{n}  \eta^{-\frac{1}{q+1}}.$$
If $\eta \leq 1 $ is   satisfied,  the solution is advanced with $u_{n+1}$ and the new step size is chosen as $h_{n+1} = h_{new}$. Often some user defined safety parameter $\beta \in (0,1) $ is introduced,
so that the condition to accept the new time step will be
$\eta \leq \beta$. Otherwise, the step is rejected and the computation is repeated with $h_{n} = h_{new}$.
In practice,
the optimal step size is calculated using appropriate safety factors,
so that
\begin{equation}
 \label{eq:stepSizeControl_04}
	h_{new} = h_{n}   \min\{\alpha_{max} , \max\{ \alpha_{min}, \alpha \eta^{-\frac{1}{q+1}}\}\}.
\end{equation}
Reasonable default values of these safety factors can be selected as $\alpha_{max} = 1.2$, $\alpha_{min} = 0.5 $, and $\alpha = 0.9$. It is important to remark that, in the application of these definitions in a multi-rate framework, the value of $\alpha_{min} $ has an impact on the maximum number of sub-steps chosen for the active components,
so that it may have to be adjusted depending on the
time scales involved in the specific problem to be solved.

In order to define the multi-rate extension of a given RK method,  we  denote by 
   $${u}_{n+1}= {\cal S} ({u}_{n}, h_n)  ,  \ \ {\cal S} :\mathbb{R}^N\times \mathbb{R}^+\rightarrow \mathbb{R}^N $$
 the global step of the basic single-rate method.
We will denote as ${\cal V}^f_n \subset \mathbb{R}^N $ the linear subspace of  fast variables   at time level $n$.
It is implicitly assumed that $d_n = \dim({\cal V}^f_n ) \ll N$. We also denote as
$P^f_n :\mathbb{R}^N\rightarrow \mathbb{R}^{d_n}$ the projector on the linear subspace of  fast variables, while
 $ {P}^{s}_n : \mathbb{R}^N\rightarrow \mathbb{R}^{N-d_n} $ denotes the projection  on the linear subspace of  slow variables $ {\cal V}_n^{s}$.
 For brevity, we will also denote ${u}^s_{n}={P}^{s}_n u_n $  and  ${u}^f_{n}={P}^{f}_n u_n$. 
The operator  
  $${u}^f_{n+1}={\cal S}^f({u}^f_{n},w, h_n )  :\mathbb{R}^{d_n}\times \mathbb{R}^{N-{d_n}}\rightarrow \mathbb{R}^{d_n} $$
  will then denote the application of the  basic method ${\cal S}  $   to  the subsystem of   $y^{\prime}=f(y,t) $ obtained  by 
  projecting ${u}_{n}$ onto ${\cal V}_n^{f} $  and assuming  that the remaining components of ${u}_{n}$
   belonging to  ${\cal V}_n^{s} $ are given by the given vector $w.$ Furthermore,
  setting  $\zeta=t_n+\tau h_{n} $ for each $\tau\in [0,1],$ we will denote as $Q(\tau) $
 a  generic  operator that interpolates known values
 of slow variables at 
  time level $\zeta. $  
{   Notice that, contrarily to what is stated in \cite{kvaerno:2000}, the use of interpolation methods at this stage does not alter the one-step nature of the resulting multi-rate methods, as will be better explained in the following.}
   This operator is used to provide intermediate values of the  slow variables
 for the application of ${\cal S}^f $ and can be given, for example, by the dense output approximation associated to a given ODE method, see e.g. \cite{hairer:2000} for relevant examples. Given these definitions,
 the general multi-rate approach we will implement can be defined as in the following
 pseudo-code.
 \vskip 0.5cm
{
\footnotesize
{\bf Multirate algorithm:}
\textit{
\begin{enumerate}
\item[1)]
Compute a tentative global step $u_{n+1}= 
{\cal S} ({u}_{n},h_n)$ and the additional approximation $\hat{u}_{n+1}$.  
\item[2)]
For each component $u^{n+1}_i, i=1,\dots,N, $ compute the error estimation $\eta^{n+1}_{i}$ as described in expression \eqref{eq:stepSizeControl_01}. 
\item[3)] Let ${\tilde\eta}^{n+1}_{i}$ 
denote the components of a vector obtained by sorting $\eta_{n+1}$ in descending order
and $s=s(i)$, $i=1,\dots,N$, the map that 
assigns to each component of $\eta_{n+1}$
the index of its location in ${\tilde\eta}_{n+1}$. Define $m \in \{1,\dots,N\} $ as the only integer such that
$m/N \leq \phi < (m+1)/N, $ where 
$\phi \in (0,1) $ is a user defined parameter 
determining the maximum fraction of fast variables that will be allowed in the simulation.
Define also the index subset
$$
{\cal S}_{n} = \{ i\in \{1,\dots,N\} :  s(i) > m \}.
$$
\item[4)] If $$\eta_s=\max_{i\in {\cal S}_{n}} \eta_i^{n+1} > \beta,$$ 
reject the global time step and recompute it using  the time step given by expression \eqref{eq:stepSizeControl_04} with respect to $\eta_s$.
\item[5)] If $$\eta_f=\max_{i\in {\cal S}_{n}^{\;c}} \eta_i^{n+1} \leq \beta,$$ 
accept the global time step and  compute the next   using the time step given by expression \eqref{eq:stepSizeControl_04} with respect to $\eta_s$. 
\item[6)] If $\eta_s \leq \beta $ and $ \eta_f > \beta, $ go multirate: 
\begin{itemize}
\item[6.1)] Define ${\cal V}^f_n$ as the subspace of $\mathbb{R}^N$ whose coordinates have indices $i$ such that $\eta_i^{n+1} > \beta$ and
$ {\cal V}_n^{s} = {{\cal V}^{f}_n}^{\,c} $ and 
${\cal S}_{n}^f$ the corresponding index subset. 
Notice that ${\cal S}_{n}^f \subset {\cal S}_{n}^{\;c}$.
\item[6.2)] 
Partition  the state space as $\mathbb{R}^N={\cal V}^f_n\oplus{\cal V}_n^{s}$, and  set  
$ {u}^s_{n+1}=P^{s}_{n}{u}_{n+1}  ;$ 
\item[6.3)]
Set ${u}^f_{n,0}={u}^f_n$. Compute  
$u^f_{n,l+1}= {\cal S}^f({u}^f_{n,l},{u}^s_{n,l}, h_n^f)$ with local time steps determined using time step control based on expression \eqref{eq:stepSizeControl_04} with respect to
 $$\hat \eta = \max_{i\in {\cal S}_{n}^f} \eta_i^{n,l+1}.$$
\item[6.4)] Proceed until the time of the global time step has been reached. Apply the same logic for step rejection and acceptance (see 4) and 5)) with respect to $\hat \eta$. Compute necessary values of the slow variables at intermediate time levels  using the interpolation operator $Q(.)$. Let $M_n$ denote the number of local steps taken.
\item[6.5)] Set  ${u}^f_{n+1}={u}^f_{n,M_n} $  and determine a new global time step based on expression \eqref{eq:stepSizeControl_04} with respect to $\eta_s$, then go to 1).
\end{itemize}
\end{enumerate}
}
}

As remarked previously, the value of $\alpha_{min}$ 
plays a role in determining the maximum number of sub-steps performed for the active variables. Therefore, it may need to be adjusted depending on specific accuracy or stability features of the method or problem under consideration. Furthermore, as discussed in Section \ref{sec:intro}, the set of active variables is chosen here in a way that has already proven to be successful in adaptive finite element techniques, see e.g. \cite{arndt:2022,
bangerth:2007, orlando:2022a, orlando:2022b}. 
This enables us to reduce or avoid conflicts
between the global time step adaptation and the multi-rate strategy, which is instead employed to achieve the same level of accuracy as the corresponding single-rate implementation at a reduced computational cost.

Last, but not least, an important remark is due for the case of strongly non-linear ODEs and implicit Runge-Kutta multi-rate integration methods.  The equations of implicit Runge-Kutta methods are usually solved using iterative Newton-Raphson methods, with an initial guess computed by extrapolation of the solution found at previous time steps. In case the RHS $f(y,t)$ of the ODEs is strongly non-linear, if the time step is too large this may cause convergence problems or even failure to converge, because the initial guess could be too far from the solution. 
This issue is much more critical for multi-rate integration, since the global time step is chosen without considering the errors of the variables belonging to the fast partition, so it is in general much larger than in the case of single-rate methods. {   For efficiency, } it is therefore essential to limit the maximum number of Newton iterations to a relatively small value such as 20, after which a shorter time step should be attempted. Otherwise, a large amount of time and computational effort could be wasted in the futile attempt to achieve convergence from an initial guess which is too far from the solution. It is also advisable to choose conservative values of $\alpha_{max}$, such as $\alpha_{max} = 1.2$, to avoid getting again into convergence issues at the next global time step.

{  
\section{Review of multi-rate stability analyses}
 \label{sec:review}
The linear stability of multi-rate methods has been studied in a number of papers. The first contribution in this respect appears to be \cite{gear:1984}. In this paper, the intrinsic difficulty of performing a stability analysis valid for the multi-rate case is acknowledged, showing that the typical reduction to a scalar problem is impossible in this case and that model problems with at least two degrees of freedom must be considered.  
A very limited stability result for the multi-rate explicit Euler method with linear interpolation is presented.

Other early attempts at similar stability analyses were presented in \cite{gunther:1993, sand:1992, skelboe:1989a, skelboe:1989b}. The results 
of these papers were reviewed and found limited and inconclusive in
\cite{kvaerno:2000}, where a first attempt at a more complete analysis was presented. In \cite{kvaerno:2000} the impact on stability of the number of sub-steps employed for the fast variables is considered. However, this is done only for a combination of implicit and explicit Euler methods for the slow and fast components, respectively.
Furthermore, in \cite{kvaerno:2000}
methods using interpolation
to provide the intermediate
values of the slow components
are excluded from consideration,
based on the incorrect consideration that these methods result in effectively two-step approaches. However, this is not the case for the methods proposed in our work, as will be clear from the analysis presented in Section
\ref{sec:stability}. Stability analyses for multi-rate versions of the BDF and linear multistep methods were then presented in
\cite{rodriguez:2004, verhoeven:2007b} and also tested on the same 2-DOF model problem.

Multi-rate methods based on implicit RK solvers and
similar in spirit to the approach pursued here were first introduced in \cite{savcenco:2007, savcenco:2008}
and their stability analysis was presented in
\cite{hundsdorfer:2009}. More specifically, in
\cite{hundsdorfer:2009} the multi-rate extension of the $\theta-$method is analyzed, employing mostly linear interpolation. However, the stability results are only obtained under rather restrictive assumptions on the general system structure and for the case of two substeps only. Furthermore,
the application to the 2-DOF model problem
presents the results as a function of parameters that are specific for the 
 $\theta-$method only. In our previous paper
 \cite{bonaventura:2020a}, we generalized the approach of \cite{hundsdorfer:2009} to the
 case of the multi-rate version of the TR-BDF2
 method, but still with the limitation to the
 two substep case.

 More recently, the multi-rate extension of 
 Generalized Additive Runge-Kutta methods
 has been proposed in \cite{gunther:2016}
 and a number of multi-rate GARK methods have been proposed in the literature, see e.g.
 \cite{gunther:2022, roberts:2020,sandu:2019,sandu:2021,sarshar:2019}.  In \cite{gunther:2016}, a very general nonlinear stability analysis was proposed for dissipative coercive problems. Furthermore, in the same work it was also shown how, conceptually, a multi-rate approach based on interpolation like that presented in this work can be reinterpreted as a variant of the more general multi-rate GARK framework. However,
 no explicit form of the amplification matrix is given for any of the methods discussed, nor
 is it explicitly shown how the interpolation procedure contributes to the amplification matrix of the resulting multi-rate method,
 thus making it difficult to apply those results in practice to the methods discussed in this paper. Therefore,  the more general nature of the analysis in \cite{gunther:2016} does not affect the novelty and the usefulness of the analysis we will present in Section
 \ref{sec:stability}.

Stability analyses for specific applications to structural mechanics problems and variational integrators were proposed, among others, in \cite{belytschko:1985,smolinski:1996, daniel:1997, daniel:1998}
and later in \cite{fong:2008, gravouil:2015}.
A nonlinear stability analyis based on the TVD property was presented in \cite{constantinescu:2007} for explicit RK methods applied to spatially partitioned discretization of hyperbolic conservation laws. Several papers have then presented applications
of multi-rate (or, in the domain specific jargon, local time-stepping) approaches to partial differential  equations, see e.g. among others
\cite{grote:2015, kang:2022, kang:2023, krivodonova:2010, schlegel:2009, seny:2014, trahan:2012}, but none of these works presents a detailed stability analysis of the methods employed.

}

\section{  Stability analysis of Runge-Kutta multi-rate methods with generic interpolation}
\label{sec:stability}

   We consider a generic linear homogeneous, constant coefficient ODE system $y^{\prime}=Ly $
  with $ y\in{\mathbb{R}}^N, $ $L\in{\cal M}_{N,N}(\mathbb{R}) $ and $\sigma(L) \subset \mathbb{R}^-\cup \{0\}.$ This assumption will allow to consider also systems with purely oscillatory behaviour, for which the application of multi-rate techniques is of great practical interest.
 We assume that the state space is partitioned \textit{a priori} as ${\mathbb{R}}^N={\cal V}_s\oplus{\cal V}_f,$
 where ${\cal V}_s $ denotes the subspace of the slow or latent variables and ${\cal V}_f $ denotes the subspace of the fast or active variables.
 We also set $d=\dim({\cal V}_f)$.
  Introducing the identity matrix ${\mathbb I}_m \in{\cal M}_{m,m}(\mathbb{R})$
 and the zero matrix ${\mathbb O}_{m,n} \in{\cal M}_{m,n}(\mathbb{R}),$ one can represent the projections onto ${\cal V}_s, {\cal V}_f $
by the matrices
 \begin{equation}
 P_s= \left[  \begin{array}{lcl}
{\mathbb I}_{N-d}   &  {\mathbb O}_{N-d,d} 
\end{array}\right] \ \ \ \ 
P_f= \left[  \begin{array}{lcl}
 {\mathbb O}_{d,N-d} &    {\mathbb I}_d
\end{array}\right] 
 \end{equation}
  and the corresponding embedding operators of ${\cal V}_s, {\cal V}_f $ into ${\mathbb{R}}^N $ by $P_s^T, P_f^T$,
  respectively, as often done in the literature on domain decomposition methods, see e.g.  \cite{quarteroni:1999}.
  As a consequence, $P_s^TP_s $ is the operator that sets to zero all the components of a vector in ${\mathbb{R}}^N $
  corresponding to the fast  variables;    $P_f^TP_f $ acts analogously on the slow variables and $P_s^TP_s +P_f^TP_f= {\mathbb I}_N$.
This entails that the model system $ y^{\prime} =L y $  can be written in terms of the partitioned matrix
\begin{equation}
  \label{eq:partsys1}
   L =\left[\begin{array}{lccl}
P_sLP_s^T& P_sLP^T_f   \\
 P_f  LP_s^T&P_f LP^T_f      
\end{array} 
\right ] =
\left[\begin{array}{lccl}
L_{ss}&L_{sf}   \\
L_{fs}&L_{ff}      
\end{array} 
\right ],
\end{equation}
or, equivalently, introducing the notation $ y_f=P_fy$, $y_s=P_sy$, as 
\begin{eqnarray}
\label{eq:partsys2}
y_s^{\prime}&=&  P_sLP_s^T y_s+  P_sLP_f^T y_f   \nonumber \\
y_f^{\prime} &=&   P_fLP_s^T y_s+  P_fLP_f^T y_f.   
\end{eqnarray}
We assume that at a given time level $t_{n} $ a time step $h=h_s$ is employed for the slow variables
and $h_f=h_s/M $ for the $l=0,\dots, M$ substeps of the fast variables, corresponding to the time levels $t_{n+l}=t_n+lh^f.$
 If  the multi-rate method is based on a generic one step time discretization method
whose amplification function can be written as
$R(z) = {N(z)}/{D(z)}$, where $N, D $ are polynomials in $z$, the value of the slow variables at the new time steps  is given by $ u_{n+1}^s=P_sR(h_sL)u_n$.
As discussed in Section \ref{sec:multirate},
values $u_{n+\tau} $ at intermediate time levels $\zeta=t_n+\tau h_s\in [t_n,t_{n+1}]$, where $\tau\in [0,1] $ denotes the
fraction of the global time step at which the interpolated value is required,  can be
computed from the values of $u_n, u_{n+1}$ by appropriate
interpolation techniques, represented formally
by a
family of operators parameterized by $\tau:$
 $$Q(\tau): \mathbb{R}^N \rightarrow\mathbb{R}^N. $$
The simplest example of such operator is given
by linear interpolation.
 Since in this case 
 $$u_{n+\tau}= (1-\tau)u_n +\tau u_{n+1} =
 (1-\tau)u_n +\tau R(h_sL) u_{n},$$ 
interpolation can be represented in this case by the
application of the operator
\begin{equation}
  \label{eq:linint} Q(\tau)  = (1-\tau) \mathbb{I}_{N} +\tau R(h_sL)
\end{equation}
to $u_n$. The slow variables at intermediate substeps are then given by $u_{n+\tau}^s= P_s Q(\tau) u_{n}$.  
 Another simple option is cubic Hermite interpolation, which is easy to implement and convenient for methods up to fourth order.
 In this case, 
also approximations of $y^{\prime}(t_{n+1}),  $ $y^{\prime}(t_n)  $  are needed, which are provided by $Lu_{n+1} $ and $ Lu_{n}$.
The interpolation operator reads in this case:
\begin{eqnarray} 
\label{eq:hermint} Q(\tau)  &=&
 (1+2\tau)  (1-\tau)^2 \mathbb{I}_{N} 
 + (3-2\tau)\tau^2 R(h_sL) \nonumber \\
 &+&h_s\tau (1-\tau)^2L
 +h_s  (\tau-1)\tau^2LR(h_sL).
\end{eqnarray}
It should be remarked, however, that
the analysis presented in \cite{kuhn:2014} shows
how Hermite interpolation may introduce
instabilities in multi-rate procedures. More accurate interpolation procedures are provided by dense output versions of the approximation methods considered, see e.g. \cite{hairer:2000} for a general discussion and
\cite{kennedy:2016} for the details of some dense output DIRK and ESDIRK methods. The
precise definition of the $Q(\tau) $ operators associated to dense output interpolators will be described in the following. {   For all these interpolation procedures,
however, as well as for the continuous output interpolators
discussed later
(see equation \eqref{opform_dense})
it should be clear that, contrarily to what stated in \cite{kvaerno:2000}, the one-step nature of the resulting multi-rate methods is preserved for RK methods.}

The goal of the stability analysis will be to represent
the discrete time evolution as $u_{n+1}= R_{mr}u_{n}, $  where $R_{mr}$ will depend on $h_s, M, $ $L, $ and on the specific properties
of the time discretization   and interpolation method employed. Since $ u_{n+1}^s=P_sR(h_sL)u_n $ and $P_s^TP_s +P_f^TP_f= {\mathbb I}_N,$
one has 
 \begin{equation}
 \label{slowvar}
 u^s_{n+1}= P_sR(h_sL)P_s^Tu^s_{n}+P_sR(h_sL)P_f^Tu^f_{n},
 \end{equation}
so that the multi-rate amplification matrix will have the form
\begin{equation}
  \label{eq:mr_amp}
 R_{mr} =
\left[\begin{array}{lccl}
P_s    R(h_sL)P_s^T  &    P_s R(h_sL)P_f^T   \\
 R_{fs} & R_{ff} 
\end{array} 
\right],
\end{equation}
where the operators $R_{fs}$ and $ R_{ff} $ depend on the specific properties of the basic single-rate method on which the
multi-rate method is based and of the interpolation operator.

 We will outline the stability analysis for generic ERK, DIRK and ESDIRK methods,
 identified as customary by their Butcher tableaux. These 
 consist of 
   $A\in {\cal M}_{s,s}(\mathbb{R})$ and    $b,c \in \mathbb{R}^s, $   where
    $s$ is the number of stages and the usual simplifying hypotheses are implicitly assumed, see e.g. \cite{hairer:2000}.
The single-rate method applied to the
problem $y^{\prime}=f(y,t)$
can therefore be written, following \cite{kennedy:2019}, as
 \begin{eqnarray}
 \label{eq:rk}
 U^{(i)}&=&u_n +h\sum_{j=1}^s a_{i,j}f(U^{(j)},t_n+c_{j}h)\ \ \ \ \ i=1,\dots,s\nonumber \\
 u_{n+1}&=&u_n +h\sum_{i=1}^sb_i f(U^{(i)},t_n+c_{i}h).
 \end{eqnarray}
In the linear case $y^{\prime}=Ly$
this yields
\begin{eqnarray}
 \label{eq:rk_linear}
 U^{(i)}&=&u_n +h\sum_{j=1}^s a_{i,j}LU^{(j)}\ \ \ \ \ i=1,\dots,s\nonumber \\
 u_{n+1}&=&u_n +h\sum_{i=1}^sb_i LU^{(i)}.
\end{eqnarray}
A dense output approximation
of $y$ at time level $\zeta=t_n+\tau h_s$, see again \cite{hairer:2000, kennedy:2016}, 
can be written as
 \begin{equation}
 \label{eq:dense1}
 u_{n+\tau}=u_n +h\sum_{i=1}^sb^*_i(\tau) LU^{(i)}, \ \ i=1,\dots,s\nonumber
 \end{equation}
 with $b^*_i(\tau)=\sum_{j=1}^{p^*}b^*_{i,j}\tau^j $ for appropriate method specific values of $p^* $
and $b^*_{i,j}$. For ERK methods
\begin{equation}
  {R}^{(1)}= \mathbb{I}  \ \ \ \ \
  {R}^{(k)}=\mathbb{I}+h\sum_{j=1}^{k-1} a_{k,j} L{R}^{(j)} , 
\end{equation}
so that $U^{(i)}={R}^{(i)}u_n$. As a consequence,
the operator representation of the dense output interpolator
can be written as
\begin{equation}
\label{opform_dense}
  Q(\tau)= \mathbb{I} +h\sum_{i=1}^sb^*_i(\tau) L{R}^{(i)}.
\end{equation}
For DIRK methods one defines instead
\begin{eqnarray}
\label{abar}
  {R}^{(1)}&=& \left(\mathbb{I}-ha_{1,1} L \right)^{-1} \nonumber \\
  {R}^{(k)}&=&\left(\mathbb{I}-ha_{k,k} L \right)^{-1}\left[\mathbb{I}+h\sum_{j=1}^{k-1} a_{k,j} L{R}^{(j)} \right ], \ \ k\geq 2,
\end{eqnarray}
so that again
the dense output interpolator can be represented by formula
\eqref{opform_dense}, where this definition of the
${R}^{(k)}$ operator is employed. By the same reasoning that allows to derive \eqref{opform_dense}, one also has that
\begin{equation}
\label{amplif_mat}
  R(hL)= \mathbb{I} +h\sum_{i=1}^sb_i L{R}^{(i)}.
\end{equation}
Notice that formulae \eqref{abar} and \eqref{amplif_mat} are actually valid also for the ERK case, by simply assuming $a_{i,i}=0$, $i=1,\dots,s$,
and in the case of ESDIRK methods, by assuming $a_{1,1}=0$. This will allow to simplify the following stability analysis,
which will be carried out explicitly for DIRK methods only, but whose results are also valid for ERK and ESDIRK methods
again with appropriate assumptions on the Butcher tableaux.

For the corresponding multi-rate methods,
denoting again by $R$ the amplification function
of the single-rate method, for the linear
problem $y^{\prime}=Ly $ the slow variables can be
computed
as discussed previously by $u^s_{n+1}= P_sR(h_sL)P_s^Tu^s_{n}+P_sR(h_sL)P_f^Tu^f_{n}$.
For the fast variables, assuming that 
 $$U^{(i,l)}_s\approx P_sQ((l+c_i)/M) u_{n} $$
for $i=1,\dots,s$,
and $l=0,\dots, M-1$,
and setting $Q^{(i,l)}=Q((l+c_i)/M)$ for brevity, one has
\begin{eqnarray}
 \label{eq:rk_mr_fast}
 U^{(i,l)}_f=u^f_{n,l} 
 &+&h_f\sum_{j=1}^s a_{i,j} L_{ff}U_f^{(j,l)}  +h_f\sum_{j=1}^s a_{i,j} L_{fs}P_sQ^{(j,l)}u_n    \nonumber \\
 u_{n,l+1}^f = u^f_{n,l} &+&h_f\sum_{i=1}^sb_i
  L_{ff} U_f^{(i,l)} + h_f\sum_{i=1}^sb_i
 L_{fs} P_s Q^{(i,l)}u_n. 
\end{eqnarray}
 
In order to derive an explicit expression for the amplification
matrix of the multi-rate method, we will now consider the specific case
of DIRK methods, which as explained before allows to derive formulae that are
also valid for ERK and ESDIRK methods.
From the basic definitions it
follows that
 
 \begin{eqnarray}
 \label{eq:dirk_mr}
 U^{(1,l)}_f&=&u^f_{n,l}+h_fa_{1,1} L_{ff}U_f^{(1,l)}+h_fa_{1,1} L_{fs}P_s
 Q^{(1,l)}u_n\nonumber \\ 
 U^{(2,l)}_f&=&u^f_{n,l} 
 +h_fa_{2,1} L_{ff}U_f^{(1,l)}+h_f a_{2,1} L_{fs}P_s
 Q^{(1,l)}u_n  \nonumber \\
 &+&h_fa_{2,2} L_{ff}U_f^{(2,l)}+h_f a_{2,2} L_{fs}P_s
 Q^{(2,l)}u_n \nonumber \\
 &\vdots& \nonumber \\
 U^{(s,l)}_f&=&u^f_{n,l} +h_f \sum_{j=1}^{s} a_{s,j} L_{ff}U_f^{(j,l)} +h_f \sum_{j=1}^{s} a_{s,j} L_{fs}P_sQ^{(j,l)}u_n \nonumber \\
 u_{n,l+1}^f&=&u^f_{n,l} +h_f
 \sum_{i=1}^sb_i  L_{ff} U_f^{(i,l)}  + h_f\sum_{i=1}^sb_i L_{fs} P_sQ^{(i,l)}u_n.   
 \end{eqnarray}
This can be rewritten as
\begin{eqnarray}
 \label{eq:dirk_mr1}
U^{(1,l)}_f&=&\left(\mathbb{I}_d-h_fa_{1,1} L_{ff} \right)^{-1} \left(u^f_{n,l} 
 +h_fa_{1,1} L_{fs}P_s
 Q^{(1,l)}u_n\right)
 \nonumber \\ 
U^{(2,l)}_f&=&\left(\mathbb{I}_d-h_fa_{2,2} L_{ff} \right)^{-1} \times \left [  u^f_{n,l}  
 + h_fa_{2,1}    \left(L_{ff}U_f^{(1,l)}+  L_{fs}P_s
 Q^{(1,l)} u_n \right )\right . \nonumber \\
  &+&  \left . \hskip 4cm h_f a_{2,2}  L_{fs}P_s
 Q^{(2,l)}u_n \right] 
 \nonumber \\
 &\vdots& \nonumber \\
 U^{(s,l)}_f&=&\left(\mathbb{I}_d-h_fa_{s,s} L_{ff} \right)^{-1}
 \times \left [
 u^f_{n,l} 
 +h_f \sum_{j=1}^{s-1} a_{s,j} L_{ff}U_f^{(j,l)} \right .\nonumber \\
 &+& \left. \hskip 4cm h_f \sum_{j=1}^{s} a_{s,j}L_{fs}P_sQ^{(j,l)}u_n 
 \right ]\nonumber \\
 u_{n,l+1}^f&=&u^f_{n,l} +h_f
 \sum_{i=1}^sb_i  L_{ff} U_f^{(i,l)}  + h_f\sum_{i=1}^sb_i L_{fs} P_sQ^{(i,l)}u_n.   
 \end{eqnarray}
We then define
\begin{equation}
  R_{ff}^{(k)}=\left(\mathbb{I}_d-h_fa_{k,k} L_{ff} \right)^{-1}\left[\mathbb{I}_d+h_f\sum_{j=1}^{k-1} a_{k,j} L_{ff}R_{ff}^{(j)} \right ]   \ \ \ \ \ k=1,\dots,s
\end{equation}
  $B^{(1,l)}=\left(\mathbb{I}_d-h_fa_{1,1} L_{ff} \right)^{-1}h_fa_{1,1}L_{fs} P_sQ^{(1,l)}   $   and for $k=2,\dots,s $
\begin{equation}
  B^{(k,l)}=\left(\mathbb{I}_d-h_fa_{k,k} L_{ff} \right)^{-1}
  \left[h_f\sum_{j=1}^{k} a_{k,j}L_{fs}P_sQ^{(j,l)}
  +h_f\sum_{j=1}^{k-1} a_{k,j} L_{ff}B^{(j,l)}\right].
  \nonumber
\end{equation}
Notice that in the ERK and ESDIRK cases this yields $B^{1,l}=\mathbb{O}_{d,N}$.
Equations \eqref{eq:dirk_mr1} can then be rewritten as
\begin{eqnarray}
 \label{eq:dirk_mr2}
 U^{(k,l)}_f&=& R_{ff}^{(k)}u^f_{n,l} + B^{(k,l)}u_n \ \ \ \ \ \ k=1,\dots,s \nonumber \\ 
 u_{n,l+1}^f&=&
 \left [ \mathbb{I}_d+h_f\sum_{i=1}^sb_i  L_{ff} R_{ff}^{(i)}\right] u^f_{n,l} \\
&+&h_f\sum_{i=1}^sb_i  L_{ff}B^{(i,l)}u_n
 + h_f\sum_{i=1}^sb_i L_{fs} P_sQ^{(i,l)}u_n  \nonumber
\end{eqnarray}
Setting then for $l=0,\dots, M-1$
\begin{eqnarray}
 C_{ff} &=&  \mathbb{I}_d+h_f\sum_{i=1}^sb_i  L_{ff} R_{ff}^{(i)} 
 \nonumber \\  
  D^{(l)} &=& \sum_{i=1}^sb_i \left [ L_{ff}B^{(i,l)}
   +   L_{fs} P_sQ^{(i,l)}\right]
\end{eqnarray}
and using the identity $P_s^TP_s +P_f^TP_f= {\mathbb I}_N$, one has for the generic step 
\begin{eqnarray}
 \label{eq:recursive}
u^f_{n+1} &=& C_{ff}  ^Mu^f_{n}  +h_f\sum_{k=1}^{M} C_{ff}  ^{M-k}
D^{(k-1)} P_f^Tu^f_{n} \nonumber \\
&+&h_f\sum_{k=1}^{M} C_{ff} ^{M-k}
D^{(k-1)} P_s^Tu^s_{n}. 
\end{eqnarray} 
As a consequence, the blocks $R_{ff}$ and $ R_{fs} $ of
the multirate amplification matrix $ R_{mr} $ in \eqref{eq:mr_amp} have the form:
\begin{eqnarray}
  \label{eq:mr_amp_rke}
  R_{ff}&=&      C_{ff} ^M + h_f\sum_{k=1}^{M} C_{ff} ^{M-k}
D^{(k-1)} P_f^T \nonumber \\
R_{fs}&=& h_f\sum_{k=1}^{M} C_{ff} ^{M-k}
D^{(k-1)} P_s^T.
\end{eqnarray}
{   Notice that, from \eqref{amplif_mat}, one has $ C_{ff} =R(h_fL_{ff})$.}

 \section{Examples of stability analysis on specific model problems}
 \label{sec:model} \indent

{   The amplification matrices derived in Section \ref{sec:stability} can be in principle computed analytically for a generic (small) linear system. Unfortunately, the resulting expressions are extremely cumbersome and not very illuminating, even in the simplest case of problems with 2 DOF.
Furthermore, using these analytic expressions, it is quite difficult to study in detail the case of large number of fast substeps, which is one of the main results of our analysis. For these reasons, we will compute the spectral radii of the stability matrices numerically.
 
In previous works devoted to the analysis of multi-rate methods
\cite{hundsdorfer:2009, kuhn:2014,kvaerno:2000, savcenco:2008},  a simple two DOF system has  been considered as model problem for multi-rate stability analysis.
This system can be written as
\begin{equation}
  \label{eq:modsys_2dof}
   L =\left[\begin{array}{cc}
-1& 1  \\
-\kappa\alpha & -\alpha 
\end{array} 
\right ]. \nonumber
\end{equation}
The parameter $\alpha>0$
determines the ratio of the
fast to slow time scale and is therefore an indicator of the stiffness of the system. The parameter $\kappa$ denotes the strength of the coupling of the fast to the slow component. For $\kappa<1$, both eigenvalues of $L$ are negative, so this is a suitable model problem for the stability analysis. 
Since it is not possible to present the results of stability analyses for a whole class of methods, we will  only discuss  some
examples of an explicit and an implicit RK method, noting that similar results have been obtained
when performing the same analysis on other methods. The main  goal  of the analysis is to assess the impact of  the multi-rate
procedure on the stability properties of the original  single-rate method. For this purpose, we will introduce
as in \cite{bonaventura:2020a} the parameter $C=h_s\Lambda, $ where $\Lambda=\max(|\lambda_i|)$ and 
$\lambda_i, i=1, 2 $ denote the eigenvalues of $L.$ This parameter, which in a PDE context would  be interpreted as an
analog of the  Courant number, must satisfy an $O(1) $ A-stability bound for single-rate explicit methods, while essentially all useful single-rate implicit methods
are A-stable for arbitrarily large values of $C.$ As a consequence, increased values of $C$ for explicit multi-rate methods will
demonstrate their potential increase in efficiency with respect to their single-rate counterparts, while  conditional stability 
for implicit multi-rate methods will demonstrate the robustness  limits of the multi-rate approach.
As example of explicit RK method, we consider the classical fourth order RK method, for which Hermite interpolation is used
in the multi-rate procedure. The Butcher tableaux of the method is reported in Appendix \ref{sec:butcher}. 
 It can be observed from Tables
 \ref{tab:stab_ERK4_alpha1_2D}-\ref{tab:stab_ERK4_alpha1000_2D} that for the explicit method considered the multi-rate approach effectively extends the stability region as the time scale separation increases. Notice that the maximum value of $C$ for which the single rate method is stable is 3.
}

\pagebreak

      \begin{table}[h!]
    \centering
    {  
    \begin{tabular}{|c|c|c|c|c|c|c|c|}
    \hline
              &$ M=2$ &  $M=4$ &         $ M=8$ &       $M=16$   &         $M=32 $     & $ M=64$      &$ M=128$ \\
             \hline
$\kappa=0.9\times10^{-5}$ & 3 &	3&         3 & 	3   &         3 &	          3   & 3 \\
\hline
$\kappa=0.9\times10^{-4}$&3 &	3 &	3& 	3& 	3  &	  3    &	3 \\
\hline
$\kappa=0.9\times10^{-3}$&3&	3 &	3 &	3 &	3&	          3 & 	3\\ 
\hline
$\kappa=0.9\times10^{-2}$ &3 &	3 &	3 &	3& 	3& 	 3  & 	3 \\
\hline
$\kappa=0.9\times10^{-1}$&4 &	4 &	4 & 	4 &         4&  	4& 	4 \\
\hline
$\kappa=0.9$ &3 &	3 &   	3& 	3& 	3& 	3& 	3\\ 
\hline
  \end{tabular}}
    \caption{{   Maximum value of $C$ for stability of the classical fourth order ERK method with Hermite interpolator,
  $\alpha=1.$}
    \label{tab:stab_ERK4_alpha1_2D}}
 \end{table}
 
         \begin{table}[h!]
 \centering
    {  
    \begin{tabular}{|c|c|c|c|c|c|c|c|}
    \hline
              &$ M=2$ &  $M=4$ &         $ M=8$ &       $M=16$   &         $M=32 $     & $ M=64$      &$ M=128$ \\
             \hline
$\kappa=0.9\times10^{-5}$ & 6 &	12&         23 & 	28   &         28 &	          28   & 28 \\
\hline
$\kappa=0.9\times10^{-4}$&6 &	12 &	23& 	28& 	28 &	  28    &	28 \\
\hline
$\kappa=0.9\times10^{-3}$&6&	12 &	23 &	27 &	26&	          26 & 	26\\ 
\hline
$\kappa=0.9\times10^{-2}$ &6 &	12 &	16 &	15& 	15& 	 15  & 	15 \\
\hline
$\kappa=0.9\times10^{-1}$&6 &	11 &	10& 	10&         10&  	10& 	10 \\
\hline
$\kappa=0.9$ &5 &	10 &   	7& 	7& 	7& 	7& 	7\\ 
\hline
  \end{tabular}}
    \caption{{  Maximum value of $C$ for stability of the classical fourth order ERK method with Hermite  interpolator,
  $\alpha=10.$}}
    \label{tab:stab_ERK4_alpha10_2D}
    
 \end{table}
         \begin{table}[h!]
 \centering
    {  
    \begin{tabular}{|c|c|c|c|c|c|c|c|}
    \hline
              &$ M=2$ &  $M=4$ &         $ M=8$ &       $M=16$   &         $M=32 $     & $ M=64$      & $ M=128$ \\
             \hline
$\kappa=0.9\times10^{-5}$ & 6 &	12&         23 & 	45   &         90 &	      $   \geq  100 $  & $\geq 100 $ \\
\hline
$\kappa=0.9\times10^{-4}$ &6 &	12 &	23& 	45& 	76  &	  74    &	74 \\
\hline
$\kappa=0.9\times10^{-3}$ &6 &	12 &	23 &	45 &	43&	          43 & 	43\\ 
\hline
$\kappa=0.9\times10^{-2}$ &6 &	12 &	23 &	25& 	25& 	 25  & 	25 \\
\hline
$\kappa=0.9\times10^{-1}$ &6 &	12 &	16 & 	15 &         15&  	15& 	15 \\
\hline
$\kappa=0.9$ &6 &	10 &   	10& 	10& 	10& 	10& 	10\\ 
\hline
  \end{tabular}}
    \caption{{  Maximum value of $C$ for stability of the classical fourth order ERK method with Hermite interpolator,
  $\alpha=100.$}}
    \label{tab:stab_ERK4_alpha100_2D}
 \end{table}

\pagebreak
 
        \begin{table}[h!]
 \centering
    {  
    \begin{tabular}{|c|c|c|c|c|c|c|c|}
    \hline
              &$ M=2$ &  $M=4$ &         $ M=8$ &       $M=16$   &         $M=32 $     & $ M=64$      &$ M=128$ \\
             \hline
$\kappa=0.9\times10^{-5}$ & 6 &	12&         23 & 	45   &         90 &	      $  \geq   100  $ & $\geq 100$ \\
\hline
$\kappa=0.9\times10^{-4}$&6 &	12 &	23& 	45& 	90  &	$  \geq 100 $   &	$\geq 100 $ \\
\hline
$\kappa=0.9\times10^{-3}$&6 &	12 &	23 &	45 &	76&	          74 & 	74\\ 
\hline
$\kappa=0.9\times10^{-2}$ &6 &	12 &	23 &	45& 	43& 	 43  & 	43 \\
\hline
$\kappa=0.9\times10^{-1}$&6 &	12 &	23 & 	25 &         25&  	25& 	25 \\
\hline
$\kappa=0.9$ &6 &	12 &   	16& 	15& 	15& 	15& 	15\\ 
\hline
  \end{tabular}}
    \caption{{  Maximum value of $C$ for stability of the classical fourth order ERK method with Hermite interpolator,
  $\alpha=1000.$}}
    \label{tab:stab_ERK4_alpha1000_2D}
 \end{table}


{  
As an example of implicit RK methods, we consider the fourth order ESDIRK method denoted as
ESDIRK4(3)6L[2]SA in \cite{kennedy:2016}, Section 7.1.1, for which the associated continuous output interpolator
was used in the multi-rate procedure. Again,
the Butcher tableaux of the method and the coefficients required to build the associated
continuous output interpolation are reported in Appendix \ref{sec:butcher}.
Repeating the previous analysis
for this method,  the stability properties are not affected by the multi-rate extension, as long as the coupling coefficient is bounded by 1.
The same seems to hold
more generally for all the unconditionally L-stable single rate methods we have studied. No results are reported in this case since the parameter
$C$ was found to be larger than the reference value 100 in all cases.
}

{   
System \eqref{eq:modsys_2dof}, however, does not have in our opinion a clear physical interpretation and does not allow to study the stability for problems with
eigenvalues located close or on the imaginary axis, which are also very important in many applications.}

We therefore consider a system with four DOF that has a clear physical
interpretation and allows to study explicitly the dependence of the multi-rate stability on the intensity of the coupling of fast and slow variables
and on the separation of their time scales. This system is an extension of a similar system considered in \cite{bonaventura:2020a}, in which however only
a partial coupling of the slow and fast variables was considered. More specifically, consider the second order system
\begin{eqnarray}
\label{eq:masses1}
m_1u_1^{\prime\prime}&=&-k_1u_1 -k_2(u_1-u_2)-c_1 u_1^{\prime}    \nonumber \\
m_2u_2^{\prime\prime} &=&  \ \ k_2(u_1-u_2)   -c_2 u_2^{\prime}.
\end{eqnarray}
It is easy to see that the system describes the dynamics of two point masses $m_1, m_2$, the first of which is subject
to elastic forces due to two springs of elastic constants $k_1, k_2$, respectively, the first of which is attached to a wall,
while the second ties the two masses. Both masses are also subject to frictional forces defined by the coefficients $c_1, c_2$.
Assuming that $y_1=u_1$, $y_2=u_1^{\prime}$, $y_3=u_2$, $y_4=u_2^{\prime}$, one can rewrite the system in first order form
as
\begin{eqnarray}
\label{eq:masses4}
y_1^{\prime}& =& \ \ y_2 \nonumber \\
y_2^{\prime}&=&-\frac{k_1}{m_1}y_1 -\frac{k_2}{m_1}(y_1-y_3)-\frac{c_1}{m_1} y_2    \nonumber \\
y_3^{\prime}& =&\ \  y_4 \nonumber \\
y_4^{\prime}&=&\ \ \frac{k_2}{m_2}(y_1 -y_3)-\frac{c_2}{m_2} y_4.
\end{eqnarray}
The system can be further rewritten to highlight the intensity of the coupling and the time scale separation between the slow and fast variables.
For $i=1,2$ we define the natural periods of the two masses as $\omega_i=\sqrt{k_i/m_i} $
and the quantities $\gamma_i=c_i/m_i$
and we set
\begin{equation}
\alpha=\frac{\omega_2}{\omega_1}  \ \ \  \ \ \beta=\frac{\gamma_2}{\gamma_1} \ \ \ \ \ \kappa=\frac{m_2}{m_1}.
\end{equation}
Notice that $\alpha, \beta $ determine the ratios of the proper periods of  the two (decoupled) masses and   of the effective frictional
forces, so that in order for $y_1,y_2$ to represent slow variables either $\alpha$, or $\beta$, or both must be significantly larger than one. On the other hand, $\kappa$ represents
the intensity of the coupling, so that for $\kappa \rightarrow 0 $ the slow dynamics tends to be decoupled from the fast dynamics,
while stronger coupling arises if $\kappa $ is of order 1. With these definitions, the system can be rewritten as
\begin{eqnarray}
\label{eq:masses5}
y_1^{\prime}& =& \ \ y_2 \nonumber \\
y_2^{\prime}&=&-\omega_1^2(1+\alpha^2\kappa) y_1 -\gamma_1 y_2+\kappa\alpha^2\omega_1^2 y_3    \nonumber \\
y_3^{\prime}& =&\ \  y_4 \nonumber \\
y_4^{\prime}&=&\ \  \alpha^2\omega_1^2y_1 - \alpha^2\omega_1^2y_3 - \beta\gamma_1 y_4.
\end{eqnarray}
In matrix form, the system can be written as $y^{\prime}
=Ly $ with system matrix given by
\begin{equation}
  \label{eq:modsys_fless}
   L =\left[\begin{array}{lccl}
L_{ss}&L_{sf}   \\
L_{fs}&L_{ff}      
\end{array} 
\right ]=\left[\begin{array}{lccccl}
0& 1& 0 & 0 \\
-\omega_1^2(1+\alpha^2\kappa) & 
-\gamma_1  &\kappa\alpha^2\omega_1^2& 0 \\
 0& 0& 0 & 1 \\
 \alpha^2\omega_1^2 & 0 & -\alpha^2\omega_1^2  &-\beta\gamma_1
\end{array} 
\right ]. \nonumber
\end{equation}
{   As for the previous model problem,
 we will perform the stability
 analysis fixing   the time scale ratios and considering
 different values of the coupling parameter $\kappa $ 
 and of the number of sub-steps $M.$} More specifically,
 we
consider time scale ratios $\alpha, \beta $ of order 10 and 100 and for each value of  $\alpha, \beta, $ we   consider
a range of values of $\kappa\in[0,1] $ and a range of values of the number of sub-steps for the fast variables. It should be remarked that the range for which the application of multi-rate methods is meaningful is $\kappa < 1$. Indeed, for values $\kappa \approx 1 $
both degrees of freedom exhibit fast dynamics, so that no actual time scale separation exists and the application of multi-rate methods would not be advantageous. Furthermore,
it is of interest to assess how the time resolution used for the fast variables affects the overall stability of the multi-rate method.

{    We consider the same examples of explicit and implicit RK methods as for the 2 DOF case. We first consider
the purely oscillatory case $\gamma_1=0$.
In this setting, the multi-rate extension does not seem able to increase substantially the stability region of the single rate scheme. Indeed, in the case of implicit methods, the multi-rate version becomes only conditionally stable with values of the
parameter $C$ close to 1.}
 
{  
If instead the real part of the system eigenvalues is non zero, the multi-rate extension of the RK method has a
significantly wider stability range. 
  As in the previous case, the results are reported in tables
  \ref{tab:stab_ERK4_gam01_alpha1}-\ref{tab:stab_ERK4_gam01_alpha1000} for $\gamma=0.01, \beta=1$ and
  $\alpha$ ranging from 1 to 1000.
It can be observed that stability is achieved for much larger values of the
parameter $C $ than those of the single-rate method,  which is approximately 3.
As expected, the stability is improved  in the weak coupling limit and if a larger number of substeps
  is employed. However, it can be observed that the maximum value of $C$ for which stability is maintained does not grow monotonically
  as a function of the number of substeps.}
  
  

  
 {   Similar results, not shown, are also obtained for the smaller value $\gamma_1=0.001 $ and for the cases 
   $\alpha=100, \beta=10  $ and $\alpha=100, \beta=100. $  
   Instead, when the time scale separation is only due
   to the real part of the eigenvalues,  such as for example in the cases
    $\gamma_1=0.01,$ $\alpha=1 $ and $\beta=10 $ or $\beta=100, $ the multi-rate extension of the explicit RK method has
   the same range of stable $C$ values as  its single-rate counterpart.}


  \begin{table}[h!]
    \centering
 {  
    \begin{tabular}{|c|c|c|c|c|c|c|c|}
    \hline
              &$ M=2$ &  $M=4$ &         $ M=8$ &       $M=16$   &         $M=32 $     & $ M=64$      &$ M=128$ \\
             \hline
$\kappa= 10^{-5}$ & 3 &	3&         3 & 	3   &         3 &	          3   & 3 \\
\hline
$\kappa=  10^{-4}$&3 &	3 &	3& 	3& 	3  &	  3    &	3 \\
\hline
$\kappa=  10^{-3}$&3&	3 &	3 &	3 &	3&	          3 & 	3\\ 
\hline
$\kappa=  10^{-2}$ &3 &	3 &	3 &	3& 	3& 	 3  & 	3 \\
\hline
$\kappa=  10^{-1}$&4 &	4 &	4 & 	4 &         4&  	4& 	4 \\
\hline
$\kappa=1$ &3 &	3 &   	3& 	3& 	3& 	3& 	3\\ 
\hline
  \end{tabular}}
    \caption{  Maximum value of $C$ for stability of the classical fourth order ERK method with Hermite  interpolator,
    $\gamma_1=0.01,$ $\alpha=1, $ $\beta=1.$}
    \label{tab:stab_ERK4_gam01_alpha1}
 \end{table}

  \begin{table}[h!]
    \centering
    {  
    \begin{tabular}{|c|c|c|c|c|c|c|c|}
    \hline
              &$ M=2$ &  $M=4$ &         $ M=8$ &       $M=16$   &         $M=32 $     & $ M=64$      &$ M=128$ \\
             \hline
$\kappa=10^{-5}$ & 6 & 12  & 23 & 	29    &        29 &	         26   & 26\\
\hline
 $\kappa=10^{-4}$ & 6 &	12 &	23 & 	14 & 
 13  &	  13    &	13 \\
\hline
$\kappa=10^{-3}$ &6 &	10 &	7 &	6 &	6&	        6 & 	6 \\ 
\hline
$\kappa=10^{-2}$ &5 &	6 &	6&	6& 	6& 	 6 & 	6 \\
\hline
$\kappa=10^{-1}$&4 &	5 &	5 & 	5 &         5&  	5& 	5 \\
\hline
$\kappa=1$ &4 &	4 &   	4& 	4& 	4& 	4& 	4\\ 
\hline
  \end{tabular}
  }
    \caption{Maximum value of $C$ for stability of the classical fourth order ERK method with Hermite  interpolator,
    $\gamma_1=0.01,$ $\alpha=10, $ $\beta=1.$}
    \label{tab:stab_ERK4_gam01_alpha10}
 \end{table}

\pagebreak


  \begin{table}[h!]
    \centering
    {  
    \begin{tabular}{|c|c|c|c|c|c|c|c|}
    \hline
              &$ M=2$ &  $M=4$ &         $ M=8$ &       $M=16$   &         $M=32 $     & $ M=64$      &$ M=128$ \\
             \hline
$\kappa=10^{-5}$ & 6 & 12  &      23 & 	
41   &       62&	          13   & 
13\\
\hline
$\kappa=10^{-4}$&6 &	12  & 22 &  7& 	6  &	  6   & 6 \\
\hline
$\kappa=10^{-3}$& 6 &	12 &	7 &	6 &	6&	         6& 	6 \\ 
\hline
$\kappa=10^{-2}$ &6 &	6 &	6 &	6& 	6& 	 6 & 	6 \\
\hline
$\kappa=10^{-1}$ &5 &	5 &	5 & 	5 &         5 &  	5& 	5 \\
\hline
$\kappa=1$ & 4 &	4 &   	4 & 	4 & 	4& 	4& 	4 \\ 
\hline
  \end{tabular}
  }
    \caption{Maximum value of $C$ for stability of the classical fourth order ERK method with Hermite  interpolator,
    $\gamma_1=0.01,$ $\alpha=100, $ $\beta=1.$}
    \label{tab:stab_ERK4_gam01_alpha100}
 \end{table}

\begin{table}[h!]
    \centering
    {  
    \begin{tabular}{|c|c|c|c|c|c|c|c|}
    \hline
              &$ M=2$ &  $M=4$ &         $ M=8$ &       $M=16$   &         $M=32 $     & $ M=64$      &$ M=128$ \\
             \hline
$\kappa=10^{-5}$ & 6 & 12  &      23 & 	
46   &       6&	          6  & 
6\\
\hline
$\kappa=10^{-4}$&6 &	12  & 23 &  6& 	6  &	  6   & 6 \\
\hline
$\kappa=10^{-3}$& 6 &	12 &	7 &	6 &	6&	         6& 	6 \\ 
\hline
$\kappa=10^{-2}$ &6 &	6 &	6 &	6& 	6& 	 6 & 	6 \\
\hline
$\kappa=10^{-1}$ &5 &	5 &	5 & 	5 &         5 &  	5& 	5 \\
\hline
$\kappa=1$ & 4 &	4 &   	4 & 	4 & 	4& 	4& 	4 \\ 
\hline
  \end{tabular}
  }
    \caption{Maximum value of $C$ for stability of the classical fourth order ERK method with Hermite  interpolator,
    $\gamma_1=0.01,$ $\alpha=1000, $ $\beta=1.$}
    \label{tab:stab_ERK4_gam01_alpha1000}
 \end{table}

{   We then consider again the ESDIRK4(3)6L[2]SA  method,   with the associated continuous output interpolator.
 As for explicit methods, if the real part of the system eigenvalues is non zero, the multi-rate extension of the ESIDRK method
 maintains unconditional stability for a wider range of values for the coupling parameter $\kappa, $
 as it can be seen from tables \ref{tab:stab_ESDIRK4_gam01_alpha1_beta1}, \ref{tab:stab_ESDIRK4_gam01_alpha100_beta1}
  for the case   $\gamma_1=0.01, \beta=1$ and $\alpha  $ ranging from 1 to 1000.
  On the other hand, values of $C=O(1) $ are   obtained for $\kappa \geq 10^{-3}$ in the
  $\alpha=10, 100  $ cases and
   for $\kappa \geq 10^{-5}$ in the $\alpha=1000  $ case. It is also to be remarked that, when unconditional stability is lost,   the number
   of sub-steps does not have a clear impact on the maximum value of $C$ for which stability is guaranteed. Similar results, not shown, are obtained with $\gamma_1=0.001,$ 
  If  instead $\gamma_1=0.01,$ $\alpha=1,  $  $\beta=10 $ or  $\gamma_1=0.01,$ $\alpha=1,  $  $\beta=100, $ so that
  time scale separation only depends on the real part of the eigenvalues, unconditional stability is maintained for
  all values of $\kappa <10^{-1}$.
}

  \begin{table}[h!]
    \centering
    {  
    \begin{tabular}{|c|c|c|c|c|c|c|c|}
    \hline
              &$ M=2$ &  $M=4$ &         $ M=8$ &       $M=16$   &         $M=32 $     & $ M=64$      &$ M=128$ \\
             \hline
$\kappa=10^{-5}$ & $\geq100 $  & $\geq100 $ &        $\geq100 $ & 	$\geq100 $  &         $\geq100 $ &	          $\geq100 $	   & $\geq100 $ \\
\hline
$\kappa=10^{-4}$& $\geq100 $  &	$\geq100 $ &	$\geq100 $ & 	$\geq100 $ & 	$\geq100 $  &	 $\geq100 $  &	$\geq100 $ \\
\hline
$\kappa=10^{-3}$& $\geq100 $  &	$\geq100 $  &	$\geq100 $ &	$\geq100 $  &	$\geq100 $ &	       $\geq100 $ & 	$\geq100 $ \\ 
\hline
$\kappa=10^{-2}$ & $\geq 100 $  &$\geq 100 $  &	$\geq 100 $  &	$\geq 100 $ & $\geq 100 $ & 	$\geq 100 $   & $\geq 100 $ \\
\hline
$\kappa=10^{-1}$& $\geq 100 $  &	7 &	7 & 	7 &         7&  	7& 	7\\
\hline
$\kappa=1$ &4 &	4 &   	      4& 	       4& 	4& 	4& 	4\\ 
\hline
  \end{tabular}
  }
    \caption{  Maximum value of $C$ for stability of fourth order ESDIRK method with continuous output interpolator,
    $\gamma_1=0.01,$ $\alpha=1,  $  $\beta=1.  $ }  \label{tab:stab_ESDIRK4_gam01_alpha1_beta1}
 \end{table}

  \begin{table}[h!]
    \centering
    \begin{tabular}{|c|c|c|c|c|c|c|c|}
    \hline
              &$ M=2$ &  $M=4$ &         $ M=8$ &       $M=16$   &         $M=32 $     & $ M=64$      &$ M=128$ \\
             \hline
$\kappa=10^{-5}$ & $\geq100 $  & $\geq100 $ &        $\geq100 $ & 	$\geq100 $  &         $\geq100 $ &	          $\geq100 $	   & $\geq100 $ \\
\hline
$\kappa=10^{-4}$& $\geq100 $  &	$\geq100 $ &	$\geq100 $ & 	$\geq100 $ & 	$\geq100 $  &	 $\geq100 $  &	$\geq100 $ \\
\hline
$\kappa=10^{-3}$& $\geq100 $  &	$\geq100 $  &	$\geq100 $ &	$\geq100 $  &	$\geq100 $ &	       $\geq100 $ & 	$\geq100 $ \\ 
\hline
$\kappa=10^{-2}$ & $\geq 100 $  &	5.0 &	5.0 &	5.0& 	5.0& 	 5.0  & 	5.0\\
\hline
$\kappa=10^{-1}$&3.0 &	3.0 &	3.0 & 	3.0 &         3.0&  	3.0& 	3.0\\
\hline
$\kappa=1$ &2.0 &	2.0 &   	       2.0& 	       2.0& 	2.0& 	2.0& 	2.0\\ 
\hline
  \end{tabular}
    \caption{Maximum value of $C$ for stability of fourth order ESDIRK method with continuous output interpolator,
    $\gamma_1=0.01,$ $\alpha=10,  $  $\beta=1.  $ }
    \label{tab:stab_ESDIRK4_gam01_alpha10_beta1}
 \end{table}

  \begin{table}[h!]
    \centering
      \begin{tabular}{|c|c|c|c|c|c|c|c|}
    \hline
              &$ M=2$ &  $M=4$ &         $ M=8$ &       $M=16$   &         $M=32 $     & $ M=64$      &$ M=128$ \\
             \hline
             $\kappa=10^{-5}$ & $\geq100 $  & $\geq100 $ &        $\geq100 $ & 	$\geq100 $  &         $\geq100 $ &	          $\geq100 $	   & $\geq100 $ \\
\hline
$\kappa=10^{-4}$& $\geq 100 $  &	$\geq 100 $ &	$\geq 100 $ & 	$\geq 100 $ & 	$\geq 100 $  &	 $\geq 100 $  &	$\geq 100 $ \\
\hline
$\kappa=10^{-3}$& $\geq 100 $ &$\geq 100 $  &	5.0 &	5.0 &	5.0&	          5.0 & 	5.0 \\ 
\hline
$\kappa=10^{-2}$ &
$\geq 100 $  &	3.0 &	3.0 &	3.0& 	3.0& 	 3.0  & 	3.0\\
\hline
$\kappa=10^{-1}$&2.0 &	2.0 &	2.0 & 	2.0 &         2.0&  	2.0& 	2.0\\
\hline
$\kappa=1$ &1.0 &	1.0 &   	       1.0& 	       1.0& 	1.0& 	1.0& 	1.0\\ 
\hline
  \end{tabular}
    \caption{Maximum value of $C$ for stability of fourth order ESDIRK method with continuous output interpolator,
    $\gamma_1=0.01,$ $\alpha=100,  $   $\beta=1.  $ }
    \label{tab:stab_ESDIRK4_gam01_alpha100_beta1}
 \end{table}

\begin{table}[h!]
    \centering
    {  
      \begin{tabular}{|c|c|c|c|c|c|c|c|}
    \hline
              &$ M=2$ &  $M=4$ &         $ M=8$ &       $M=16$   &         $M=32 $     & $ M=64$      &$ M=128$ \\
             \hline
             $\kappa=10^{-5}$ & $\geq100 $  & $\geq100 $ &        $\geq100 $ & 	$\geq100 $  &         $\geq100 $ &	          $\geq100 $	   & $\geq100 $ \\
\hline
$\kappa=10^{-4}$& $\geq 100 $  &	$\geq 100 $ &	5 & 	5 & 	 5  &	
5  &	5 \\
\hline
$\kappa=10^{-3}$& 3 &3  &	3 &	3 &	3&	          3 & 	3 \\ 
\hline
$\kappa=10^{-2}$ &
2 &	2 &	2 &	2& 	2& 	 2  & 	2\\
\hline
$\kappa=10^{-1}$&1 &1 &	1 & 	1 &     1&  	1& 	1\\
\hline
$\kappa=1$ &1 &	1 &   	       1& 	       1& 	1& 	1& 	1\\ 
\hline
  \end{tabular}
  }
    \caption{Maximum value of $C$ for stability of fourth order ESDIRK method with continuous output interpolator,
    $\gamma_1=0.01,$ $\alpha=1000,  $   $\beta=1.  $ }
    \label{tab:stab_ESDIRK4_gam01_alpha1000_beta1}
 \end{table}

\begin{figure}
    \centering
    \includegraphics[width=0.45\linewidth]{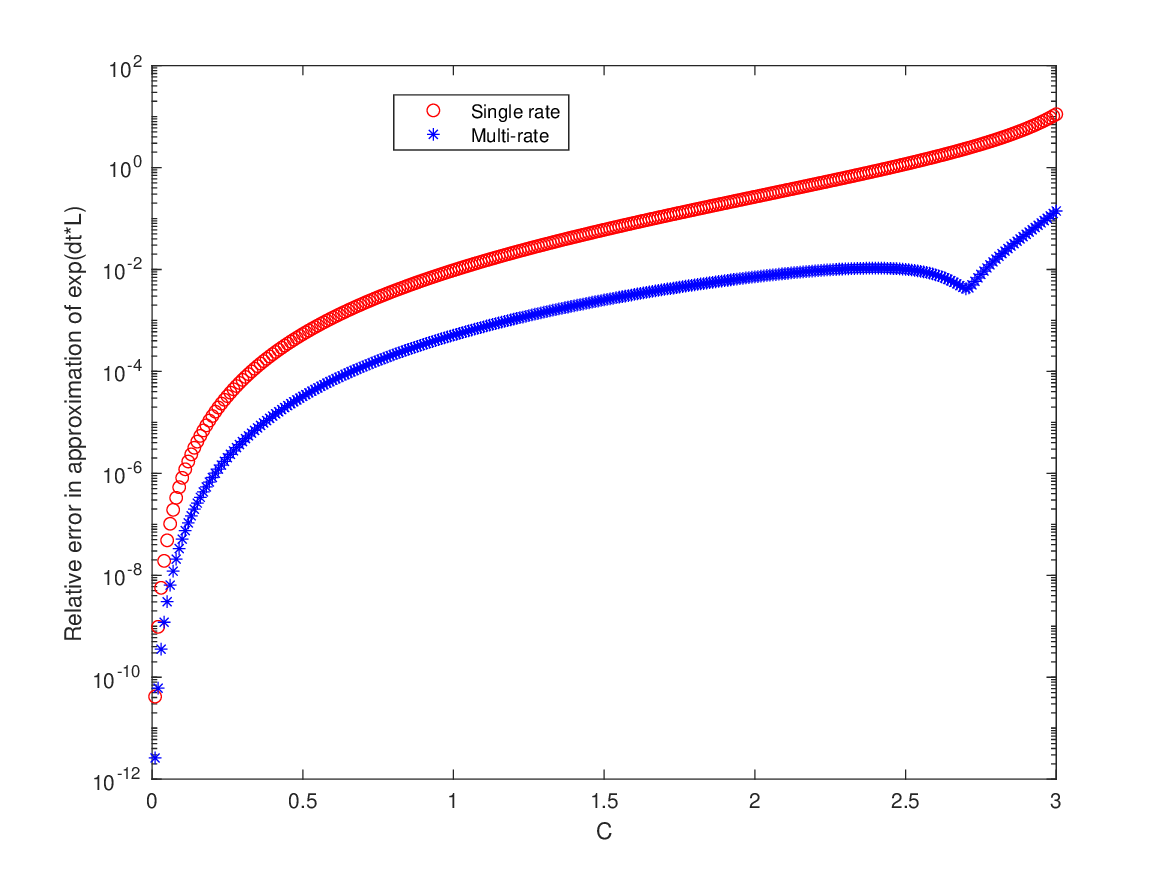}a)
    \includegraphics[width=0.45\linewidth]{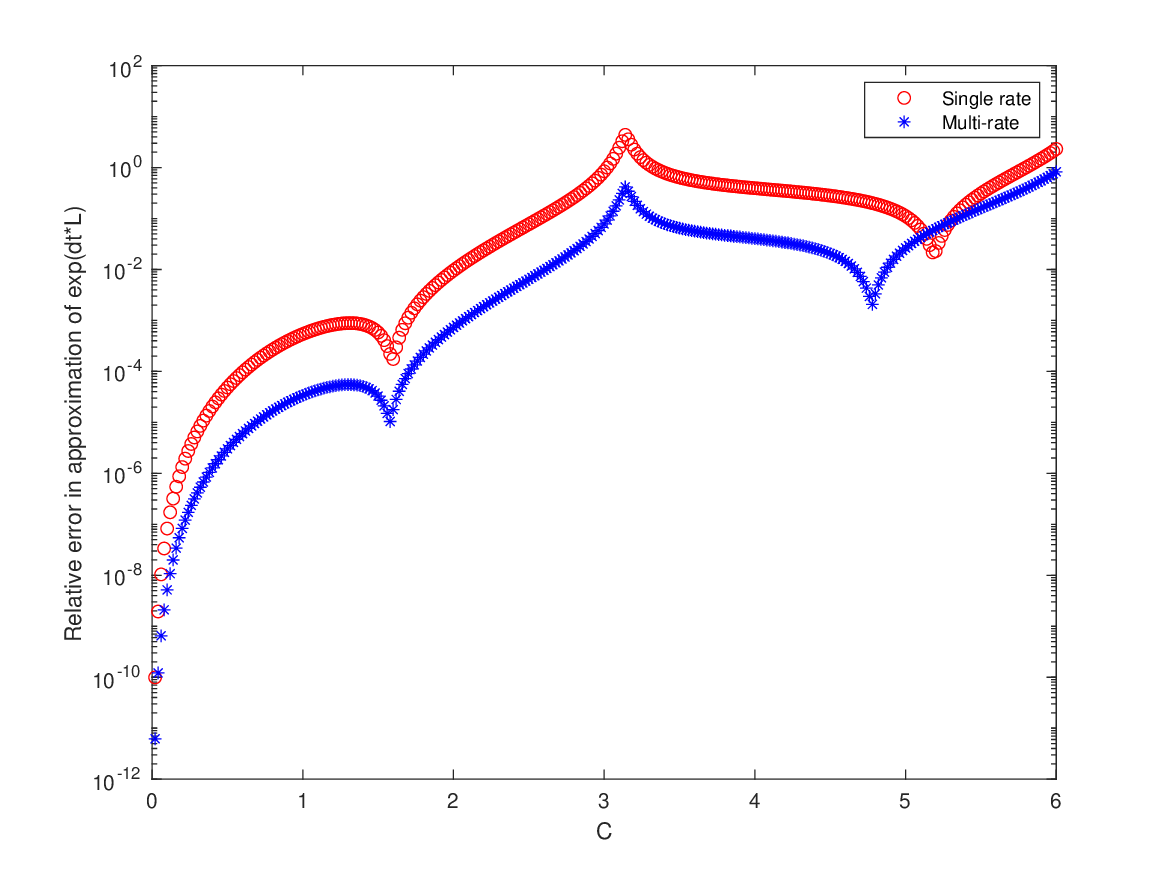}b)
    \caption{Errors with respect of exact solution at increasing values of $C$ for a) explicit  fourth order RK method, b) ESDIRK fourth order method.}
    \label{fig:accuracy}
 \end{figure}

The previously introduced expression for the multi-rate amplification matrix also allows to assess the accuracy gain achieved by the use
of smaller time steps on the fast components. For example, we consider the case of a
system with $\gamma= 0.01, $ $\alpha=50, \beta=1$ and weak coupling $\kappa=10^{-3}$. 
We compute the exact evolution matrix $\exp(Lt) $ of problem $y^{\prime}=Ly $
and compare it to appropriate
powers of the single-rate and multi-rate amplification matrix for different values of the time step. For the multi-rate methods, $M=10$ substeps were employed.
The relative errors in the $l^2$ operator
norm  of the
approximations of   $\exp(Lt) $ obtained by the  explicit fourth order RK method and fourth order ESDIRK method considered
before are reported in Figure \ref{fig:accuracy}, for values of the time step corresponding to increasing values of the parameter $C$.
It can be observed that significant improvements in the accuracy of the approximation obtained with a given time step can be achieved.

{  
 \section{Implementation in the OpenModelica  environment}
\label{sec:openmod} \indent

The multi-rate Runge-Kutta methods presented in the previous Sections were implemented in the GBODE solver, which is part of the open-source OpenModelica simulation environment described in  \cite{fritzson:2020}.

OpenModelica allows to model and simulate dynamical systems described by differential-algebraic equations using the high-level, equation-based, object-oriented modelling language Modelica \cite{MattssonElmqvistOtterCEP1998}. 
OpenModelica applies structural analysis, symbolic manipulation and numerical solvers to reduce the (possibly high-index) differential-algebraic formulation of the system model to a set of explicit ordinary differential equations, producing efficient C code to compute the right-hand-side of Equation \eqref{eq:partitioned_sys} and its Jacobian $\partial f / \partial y$. This code is then linked to ODE solvers such as GBODE or state-of-the-art solvers such as IDA or CVODE from SUNDIALS \cite{hindmarsh2005sundials}. In order to provide an efficient implementation of the multi-rate solver, the generated C code was augmented with a selection algorithm that, based on structural dependencies, only runs the parts of the code that are strictly necessary to compute the derivatives of the fast states $f_f(y_s, y_f,t)$ of Equation \eqref{eq:partitioned_sys}, as well as the Jacobian $\partial f_f / \partial y_f$, for the computation of local refinement steps.

The implementation of the multi-rate method in GBODE is quite efficient from a computational point of view: it is written using the C language, it uses the KINSOL solver of SUNDIALS to solve the nonlinear equations of implicit methods such as ESDIRK, relying on sparse linear algebra methods. It also supports dense output to provide the solution on a regular time grid, as well as to compute accurate high-order interpolation of slow variables during refinement steps.
}

\section{Numerical results}
\label{sec:tests} \indent

In this Section we will present three numerical test cases, reporting results with the multi-rate versions of the methods fully described in Appendix \ref{sec:butcher}. The numerical results demonstrate how the proposed multi-rate method works when simulating non-trivial systems. Although OpenModelica is capable of handling extremely complex dynamical system models, for this paper we have selected test cases whose model equations can be written explicitly in compact form, so that the
results are in principle reproducible with other implementations of the multi-rate method or different numerical approaches. {  All tests were carried out on a laptop with an i7-1365U CPU running Windows 11.}
   
\subsection{Inverter Chain}
\label{sec:inverter}
The first test case consists of the model of an inverter chain, which is an important test problem for electrical circuits that has already been considered in the literature on multi-rate methods, e.g., in \cite{savcenco:2007,verhoeven:2007}. The system of equations is given by
\begin{align}
    & y'_1(t) = U_{op} - y_1(t) - \Gamma g(u(t), y_1(t)) \\
    & y'_j(t) = U_{op} - y_j(t) - \Gamma g(y_{j-1}(t), y_j(t)), \quad j = 2, 3, ...,N,
\end{align}
where $y_j$ is the output voltage of the $j$-th inverter, $u(t)$ is the input voltage of the first inverter, $U_{op}$ is the operating voltage corresponding to the logical value 1, and $\Gamma$ is a stiffness parameter. The function $g$ is defined as
\begin{equation}
    g(y, z) = \max(y - U_\tau, 0)^2 - \max(y - z - U_\tau, 0)^2,
\end{equation}
where $U_\tau$ is a switching threshold. The system has a stable equilibrium if the input is zero, the odd-numbered inverters have output zero and the even-numbered inverters have output one. If the input is increased up to $U_{op}$, the first inverter output switches from zero to one, triggering a switching cascade that propagates through the entire chain at finite speed. If the first inverter input is then switched back to zero, a second switching cascade propagates again through the system, bringing it back to the original equilibrium. This problem is obviously a good candidate for adaptive multi-rate integration, since only a small fraction of inverter outputs are changing at any given time. It is also representative of the simulation of real-life electronic digital circuit models, where only a small fraction of gates are active at any given point in time, while most other gates sit idle in a stable state.

The model was set up with $N = 1000$ inverters, $U_{op} = 5$, $U_\tau = 1$, and $\Gamma = 500$. The integration interval is $t \in [0, 200]$, with initial conditions $y_j(0) = 6.247\times 10^{-3}$ for even $j$ and $y_j(0) = 1$ for odd $j$; $u(t)$ is a continuous piecewise-linear function which is zero until $t = 5$, then increases to 5 until $t = 10$, remains constant until $t = 15$, decreases to zero until $t = 20$ and then remains zero for $t > 20$.
We have used the third order ESDIRK method denoted as ESDIRK3(2)4L[2]SA  in  \cite{kennedy:2016} to compute a reference single-rate solution with relative and absolute tolerance $10^{-9}, $ a single-rate solution with relative and absolute tolerance $10^{-5}$ and multi-rate solution with relative and absolute tolerance $10^{-5}$.
For the multi-rate algorithms, the parameter values $\alpha_{max} = 1.2$, $\alpha_{min} = 0.5$, $\alpha = 0.9$, $\phi = 0.05$, and $\beta = 1$ were used; the associated  
dense output reconstruction was used for the interpolation operator $Q(\tau)$, both in the multi-rate procedure to interpolate the slow variables when integrating the fast variables, and when re-sampling the solution over a regular time grid with 0.01 s intervals.

\begin{figure}[tbph]
    \centering
    \includegraphics{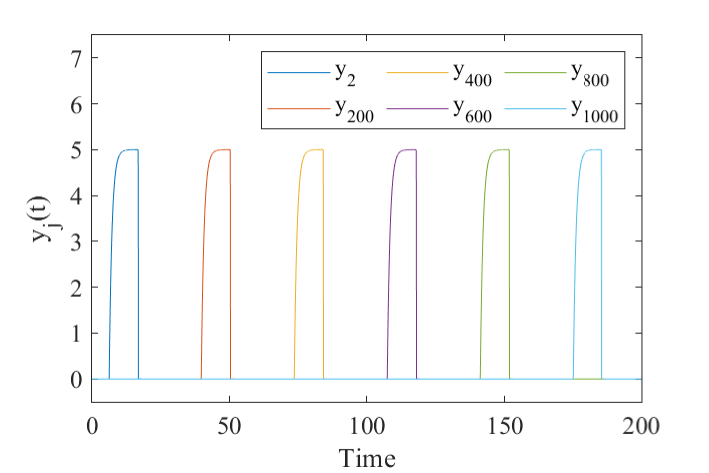}
    \caption{Reference solution for selected even inverter outputs.}
    \label{fig:inverters-reference}
\end{figure}
\begin{figure}[tbph]
    \centering
    \includegraphics{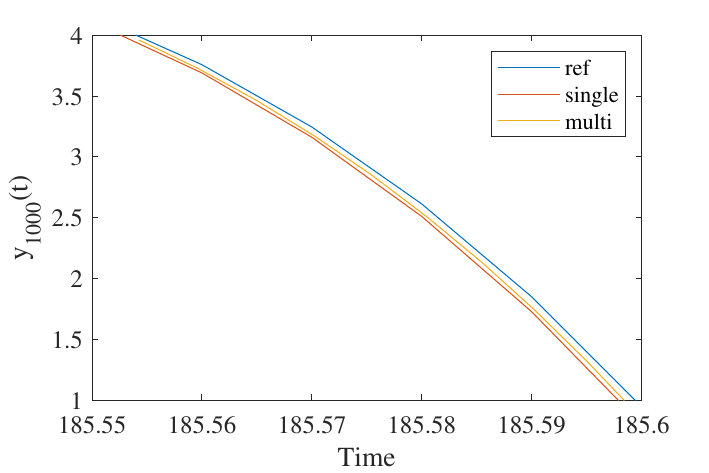}
    \caption{Comparison of reference, single-rate, and multi-rate ESDIRK3(2)4L[2]SA solutions in the inverter test case.}
    \label{fig:inverters-accuracy}
\end{figure}

Figure \ref{fig:inverters-reference} shows the output of the even inverters $y_{200}$, $ y_{400}$, $y_{600}$, $y_{800}$, $y_{1000}$ in the reference solution, which follow the rising and falling edges of the first even inverter output $y_2$ with increasing delay, corresponding to about 150 time units for the last inverter output $y_{1000}$. {   Figure \ref{fig:inverters-accuracy} shows a detail of the falling edge of the last inverter output $y_{1000}$, comparing the reference, single-rate, and multi-rate solutions. It is apparent how the multi-rate and single-rate solutions are very close to each other and that both are shifted in time with respect to the reference one by about 0.0015 time units, corresponding to a relative error of aboutn $10^{-5}$, which is compatible with the tolerances employed.}

\begin{figure}[tbph]
    \centering
    \includegraphics[width=\textwidth]{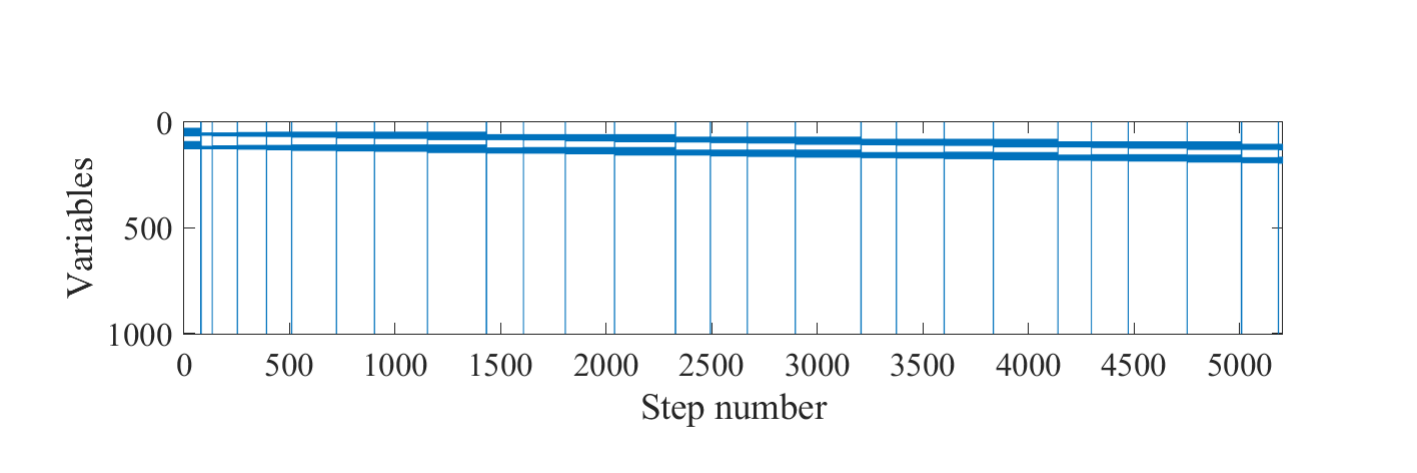}
    \caption{Activity diagram of the ESDIRK3(2)4L[2]SA multi-rate solution in the inverter test case.}
    \label{fig:inverters-activity}
\end{figure}

{   Figure \ref{fig:inverters-activity} } shows the activity diagram of the multi-rate algorithm between $t = 25$ and $t = 35$. Here, the rows correspond to the variables in vector $y$, while the columns correspond to time steps; blue points correspond to variables that are active in a given time step. The full blue vertical lines correspond to global time steps, where all variables are involved in the computation of the time step, while the slightly slanted pairs of horizontal blue segments correspond to the rising and falling inverter output waves, slowly moving along the inverter chain. Only the inverters currently undergoing the transition are automatically selected as active variables during refinement steps. The activity plot clearly highlights the ability of the self-adjusting multi-rate approach to identify automatically the set of fast evolving variables, which is essential if this set is time-varying.

In this specific example, the activity of fast states is quite localized: the typical size of the fast variable subspace is about 5, i.e., $1/200$ of the global system size $N$. This means that, as long as $\phi$ is significantly less than one (otherwise most steps would be chosen as global ones), and larger than $1/200$, the same number of fast states will be selected. This was experimentally confirmed by running the simulations with different values of $\phi$ in the range $0.4<\phi<0.01$, which led to very similar results in terms of number of global and refinement steps. In general, the performance of the proposed multi-rate algorithm will not be too sensitive to the specific value of $\phi$ for those cases where the multi-rate strategy is more advantageous, i.e., when there is strongly localized activity most of the time, so that only a tiny fraction of variables shows significant integration errors for most of the time steps.

{   Experimenting with this test case showed an important fact, that can have a substantial impact on the solver performance. The adoption of a multi-rate method on a strongly nonlinear model such as the one considered in this example has consequences on the convergence of Newton's method when solving the implicit equation for global steps. When single-rate integration is performed, the length of the steps is severely limited by the accuracy requirement; as a consequence, the extrapolation of the previous solution provides a good initial guess for Newton's method, that converges in a limited number of iterations in most cases, with a few convergence failures. When multi-rate integration is performed, global steps are much longer (100 times more on average, in this case), so the convergence of Newton's method in such steps can be more problematic. }

\begin{table}
    \centering
    {\small
      
    \begin{tabular}{|l|c|c|c|c|}
    \hline
       & Single-Rate & Multi-Rate & Multi-Rate             & IDA\\
       &             & Jac A      & Jac B                  &    \\
       \hline
       Accepted global  steps                     & 65316       &  495    &  510    &    55064 \\
       \hline
       Rejected global steps (error test failure) & 12081       &         &  0      &     8468 \\
       \hline
       Rejected global steps (convergence failure) & 1          &  103    &  107    &       18 \\
       \hline
       Accepted fast  steps                        &    //      &  66915  &  66954  &    //    \\
       \hline
       Rejected fast  steps (error test failure)   &    //      &  11654  &  11641  &   //     \\
       \hline
       Global RHS calls                            & 1221140    &  593591 &  261995 &  104003 \\
       \hline
       Global Jacobian computations                & 6943       &  24339  &  18266  &   15741 \\
       \hline
       Local RHS calls                             & //         & 1143414 &  566482 &   //    \\
       \hline
       Local Jacobian computations                 & //         &  12889  & 235762  &   //    \\
       \hline
       Simulation time in s                        & 66.1       &  37.8   & 14.6    &   12.5   \\
       \hline
    \end{tabular}
    }
    \caption{Performance of single- and multi-rate ESDIRK3(2)4L[2]SA methods and of IDA solver in the inverter test case.}
    \label{tab:inverters}
\end{table}

{    The comparison of some performance indicators of the single and multi-rate algorithms is shown in Table \ref{tab:inverters}. The first column refers to the single-rate method, while the second and the third refer to two variants of the multi-rate method. The Jac A variant uses the same Jacobian recomputation strategy of the single-rate method, i.e., the Jacobian is only recomputed by KINSOL according to its heuristic criteria, that try to strike a balance between the added computational effort of re-computing and re-factorizing the Jacobian and the added computational effort of additional iterations due to the use of outdated Jacobians. As a result, the Jacobian is pre-emptively re-computed by GBODE every 10 integration steps. In the Jac B variant, instead, GBODE pre-emptively re-computes the Jacobian at the beginning of each new time step.

Even though the total number of steps of the single-rate and multi-rate method is comparable, the multi-rate method requires less than $1\%$ global steps than the single-rate method, while $99\%$ of the steps are local refinement steps, involving only a small number of variables; thus, there is a potential for a speed-up factor up to 100.

Unfortunately, the considered test case is quite stiff and strongly non-linear, so the aforementioned issue with the convergence of the implicit stage computations is particularly severe. Whereas the average number of calls to the right hand side of Eq. \eqref{eq:partitioned_sys} per step is about 20 for the single-rate method, which is reasonable for a third-order implicit single-step method, it turns out to be over 1000 for the global steps of the multi-rate method, due to a massively increased number of iterations per step to achieve successful convergence. As a consequence, the speed-up factor turns out to be less than 2. In this case, adopting strategy Jac B is very beneficial, as it leads to a speed-up factor of about 4, due to the more than halved number of global and local right-hand-side computations that stems from the reduced number of iterations, as well as to a significant reduction of Jacobian computations for the global steps. The increased number of Jacobian computations for the local steps is probably not too relevant, due to their very small size.

The last column of Table \ref{tab:inverters} allows to compare the performance of the state-of-the art IDA solver to single-rate and multi-rate ESIDRK3 as implemented in GBODE. IDA turns out to be significantly faster than single-rate ESDIRK and slightly faster than the multi-rate ESDIRK for two main reasons. The first is that IDA is a multi-step method that requires to solve only one implicit system of equations per step, whereas ESDIRK3 is a three-stage Runge-Kutta method that requires to solve an implicit system three times per step. The second is that IDA is a highly optimized state-of-the-art solver, whereas GBODE is not yet fully optimized from the point of view of Jacobian recomputation strategies. Such an optimization, though, goes beyond the scope of the present paper.
}

\subsection{Finite difference discretization of the Burgers equation}
\label{sec:burgers}

The second  test case consists of an application to a PDE problem.  We consider the viscous Burgers equation

\begin{equation}
\frac{\partial u }{\partial t}  
+ u\frac{\partial u }{\partial x} =\nu \frac{\partial^2 u }{\partial x^2} 
\label{eq:burgers}
\end{equation}
as a simple model of computational fluid dynamics applications. The equation is considered
on the domain $[0,25] $ and time interval $[0,5].$ 
We assume $\nu=10^{-2}$  and an initial datum given by
\begin{equation}
u_0(x)=\exp\left\{-\left(\frac{x-L/2}{L/50}\right)^2\right\}.
\end{equation}
Equation \eqref{eq:burgers} is semi-discretized in space by a standard centered finite difference  approximation on a uniform mesh on $[0,25] $ with $N=1000 $ nodes. We used the third order ESDIRK method denoted as
ESDIRK3(2)4L[2]SA  in  \cite{kennedy:2016} to compute
a reference single-rate solution with relative and absolute tolerance {   $10^{-9}$}, a single-rate solution with relative and absolute tolerance $10^{-5}$, two multi-rate solutions with relative and absolute tolerance $10^{-5}$, one with $\phi = 0.2$ and one with $\phi = 0.04$, and another two multi-rate solutions with relative and absolute tolerance $10^{-6}$, one with $\phi = 0.2$ and one with $\phi = 0.04$. The associated dense output reconstruction was used for the interpolation operator $Q(\tau)$, both in the multi-rate procedure to interpolate the slow variables when integrating the fast variables, and when re-sampling the solution over a regular time grid with 0.1 s time intervals. For the multi-rate algorithms, the parameter values $\alpha_{max} = 1.2$, $ \alpha_{min} = 0.5$, $\alpha = 0.9$, and $\beta = 1$ were used. 

\begin{figure}[tbph]
    \centering
    \includegraphics{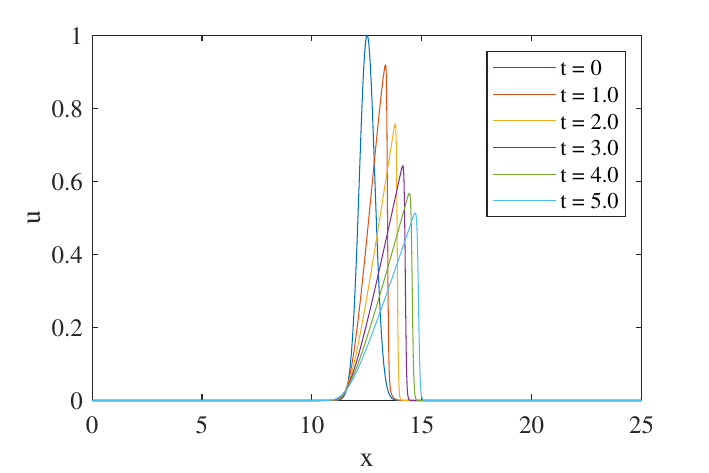}
    \caption{Reference solution for Burgers equation test case.}
    \label{fig:burgers-reference}
\end{figure}

Fig. \ref{fig:burgers-reference} shows the reference solution $u(x,t)$ along the spatial axis at five subsequent time instants. As time progresses, a sharp shock wave is formed on the right boundary of the active region, only slightly smoothed by the diffusion term, moving to the right at finite speed, followed by a trailing region where the values of $u$ slowly get back to zero. After $t = 1$, for each time instant the variable vector can be partitioned in three subsets. One, involving a few percent of the node variables, corresponds to  the fast-changing shock wave. The second, involving about 15\% of the node variables, corresponds to the slow trailing wave, while the remaining part of the node vector contains values that always remain very close to zero.

\begin{figure}[tbph]
    \centering
    \includegraphics[width=0.49\textwidth]{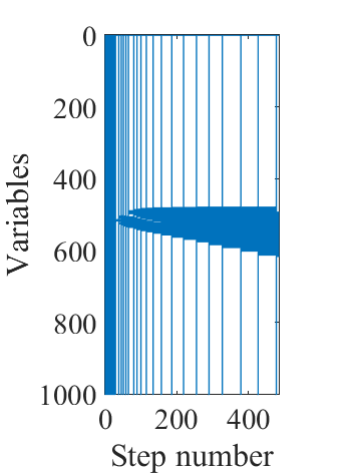}
    \includegraphics[width=0.49\textwidth]{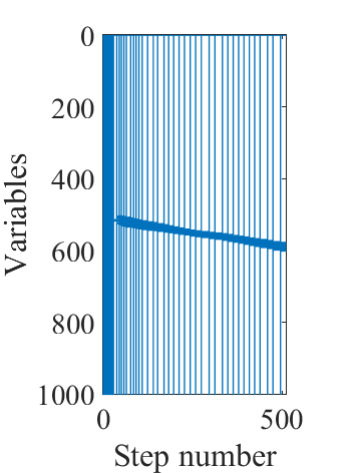}
    \caption{Activity diagram of Burgers' equation test case, tolerance $10^{-5}$, left $\phi = 0.2$, right $\phi=0.04$.}
    \label{fig:burgers-activity_1e_5}
\end{figure}

 The activity diagram in the case $\phi = 0.2$ is shown in the {   left half of
 Fig. \ref{fig:burgers-activity_1e_5}}. For this computation,  relative and absolute tolerance values were set to $10^{-5}$. {   The first 30 steps are all global steps, where the solver increases the step size from the initial value of $1.4\cdot 10^{-4}$ to $2.3\cdot 10^{-2}$, when the error tolerance threshold is first violated. Then, multi-rate integration kicks in, with}  the fast variable subsets ${\cal V}^f_n$ encompassing the first two above-mentioned subsets, where basically all the action takes place in the solution, corresponding to the central blue part of the diagram. {   After $t = 1$, which corresponds to the $130^\mathrm{th}$ step, only a very few global steps, about one every 30 fast refinement steps, are taken, essentially to capture the progressive enlarging of the fast variable subsets ${\cal V}^f_n$ as the shock-wave moves to the right of the spatial domain.}

If a smaller value $\phi = 0.04$ is chosen, the activity diagram shown in the {   right half of Fig. \ref{fig:burgers-activity_1e_5} is obtained. In this case, the sets ${\cal V}^f_n$ cannot be large enough to encompass the whole region around the shock with non constant values. Hence, they can only} cover the much narrower subset of variables corresponding to the shock transition, moving from left to right in the spatial domain, while the degrees of freedom corresponding to the trailing region and to the essentially constant part of the solution {   both end up in the subset of slow variables ${\cal V}^s_n$}. Therefore, compared to the previous solution, there are more global steps, but on the other hand the fast steps involve a much smaller set of variables.

\begin{figure}[tbph]
    \centering
    \includegraphics[width=\textwidth]{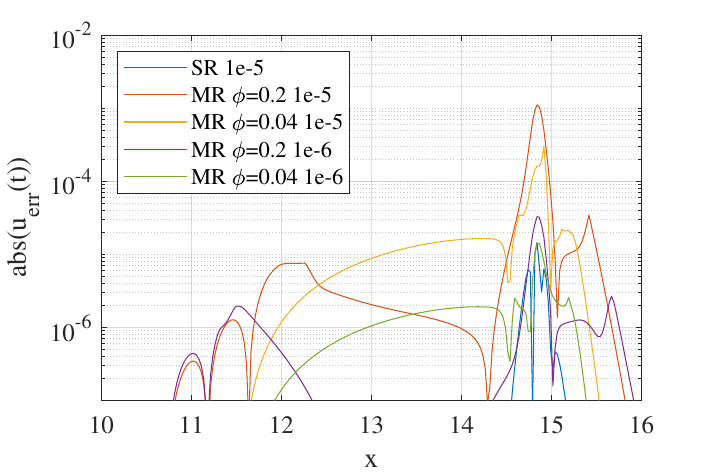}
    \caption{Absolute error for the solution of Burgers equation.}
    \label{fig:burgers-error}
\end{figure}

{   Fig. \ref{fig:burgers-error} } shows the spatial distribution of the absolute error over the spatial domain at $t = 5$. The error is computed considering the difference between the single-rate and multi-rate solutions and the reference solution computed with the single rate method employed with tolerance {  $10^{-8}$}. It is to be remarked that, in this way,  an estimate of the time discretization error is only  obtained, since all methods use the same spatial discretization. The interaction of time and space discretization errors is an extremely important point for application to PDE solvers, which is however beyond the scope of this work.   For the single-rate solution with tolerance $10^{-5}$, the maximum error is about $1.5 \cdot 10^{-5}$, which is consistent with the selected tolerance. The multi-rate solutions with tolerance $10^{-5}$ show a significantly larger error around $x = 14.8$, where the steepest region of the shock wave is located, with maximum errors of about $3 \cdot 10^{-4}$ for $\phi = 0.04$ and $10^{-3}$ for $\phi = 0.02$. {   On the other hand, it is important to consider that the error tolerance of $10^{-5}$ concerns the local error over one integration step, whereas the errors shown in Fig. \ref{fig:burgers-error} are global errors accumulated over the entire integration time interval. Most importantly, when reducing the tolerance to $10^{-6}$, the maximum absolute error is dramatically reduced to about $10^{-5}$, thus validating the correct behaviour of the multi-rate error control mechanism.}

Finally, {   Table \ref{tab:burgers} compares the performance of the single- and multi-rate simulations. As in the previous case, the overall number of steps in the multi-rate case (485) is only marginally larger than the number of steps in the single rate case (383), while the number of global steps is one order of magnitude smaller. In principle, a significant performance improvement would be expected. However, the observed simulation speed-up factor was
limited, between 1 and 2 depending on the specific settings. Note that in this case it was not necessary to pre-empt the computation of the Jacobian for each global time step to achieve the best simulation time - in fact, such a choice significantly increased it. This behaviour is likely due to the test case being markedly less nonlinear than the previous one, as also clearly indicated by the fact that there were zero convergence failures reported, so that there is a clear advantage at re-using outdated Jacobians for several global time steps.

The reason for the disappoining speed-up of the multi-rate algorithm can probably be found in the overhead of the multi-rate code implementation, which requires to set up the solver and determine the relevant section of RHS computation code to run for each fast refinement step; given the very simple expressions of the RHS terms, which require very limited computation effort, this overhead likely plays a major role in this case. The comparison of the single-rate ESDIRK3 performance with the IDA performance shows a ratio of 4.5, which can be largely attributed, as in the previous test case, to the larger number of nonlinear systems to be solved at each step (3 for ESIDIRK3 vs. 1 for IDA), as well as to the better optimizations of a state-of-the-art solver like IDA.
Further optimizations of the GBODE code will be necessary, which are, however, beyond the scope of this paper. }

\begin{table}
    \centering
    {\small
      
    \begin{tabular}{|l|c|c|c|c|c|c|}
       \hline
                                                    & SR          & MR              & MR             & MR              & MR            & IDA    \\
                                                    &             &  $\phi=0.2$     &  $\phi=0.2$   &  $\phi = 0.04$  & $\phi = 0.04$  &   //   \\ 
                                                    &             &  $tol=10^{-5}$  & $tol=10^{-6}$  &  $tol=10^{-5}$  & $tol=10^{-6}$ &  $tol=10^{-5}$    \\
       \hline
       Accepted global  steps                       & 383        &  49             &  56            & 67              & 112            &  393   \\
       \hline
       Rejected global  steps                       & 0          &  0              &  0             & 0               & 0              &  0     \\
       \hline
       Accepted fast  steps                         &    //      &  436            &  873           & 448             & 1110           &  //    \\
       \hline
       Rejected fast  steps                         &    //      &  0              &  0             & 0               & 0              &  //    \\
       \hline
       Global RHS calls                   &    5667    &  1250           &  1809          & 1579            & 3183           & 431    \\
       \hline
       Global Jacobian comp.       &    3       &  35             &  28            & 35              & 29             & 41     \\
       \hline
       Local RHS calls                    &    //      &  5390           &  8169          & 5041            & 9157           &   //    \\
       \hline
       Local Jacobian comp.        &    //      &  43             &  50            & 79              & 162            &  //     \\
       \hline
       Simulation time in s                         &    0.168   &  0.084          &  0.125         & 0.096           & 0.188          &  0.037 \\
       \hline
    \end{tabular}
    }
    \caption{Performance of single- and multi-rate ESDIRK3(2)4L[2]SA methods on Burgers' equation test case.}
    \label{tab:burgers}
\end{table}

\subsection{Thermal model of a large building heating system}
\label{sec:heating}
The third test case consists of an idealized model of the thermal behaviour of a large building with $N = 100$ heated units. For simplicity, it is assumed that each unit, with temperature $T_{u,j},$ is well-insulated from the others and only interacts with a central heating system with supply temperature $T_s$ and with the external environment at temperature $T_e$, which is assumed to vary sinusoidally during the day with a peak at 14:00. The thermal supply system has a temperature $T_s$, which is regulated by a proportional controller with fixed set point $T_s^0$, that acts on the heat input $Q_s$. The total energy consumption of the thermal supply system is tracked by the variable $E$.

Each unit has a proportional temperature controller with variable set point $T_{u,j}$, acting on the command signal $u_j$ of a heating unit with thermal conductance $G_h$ (e.g., a fan coil with variable fan speed), which responds to the command as a first-order linear system with time constant $t_h$. The heating unit then enables the heat transfer between the heating fluid of the supply system and the room, with thermal power $Q_{h,j}$. Each unit also exchanges a thermal power $Q_{e,j}$ with the external environment.

The model is described by the following system of differential-algebraic equations:
\begin{align}
& T_e = 278.15 + 8 \cos \left( 2\pi\frac{t - 14 \times 3600}{24 \times 3600} \right) &\\
& Q_s = \operatorname{sat}(K_{ps}Q_{max}(T_s^0 - T_s), 0, Q_{max}) &\\
& Q_{ht} = \sum_1^N{Q_{h,j}} & \\
& C_s \frac{dT_s}{dt} = Q_s - Q_{ht} & \\
& t_h\frac{dG_{h,j}}{dt} = u_j G_{hn} - G_{h,j}, & \quad j = 1, ..., N \\
& Q_{h,j} = G_{h,j}(T_s - T_{u,j}), & \quad j = 1, ..., N \\
& Q_{e,j} = G_u (T_{u,j} - T_e), & \quad j = 1, ..., N \\
& C_{u,j} \frac{dT_{u,j}}{dt} = Q_{h,j} - Q_{e,j}, & \quad j = 1, ..., N \\
& u_j = \operatorname{sat}(K_{pu}(T^0_{u,j} - T_{u,j}), 0, 1), & \quad j = 1, ..., N \\
& T^0_{u,j} = f_j(t), & \quad j = 1, ..., N \\
& \frac{dE}{dt} = Q_s
\end{align}
which can be easily reformulated by substitution as a system of $2N+2$ {  explicit} ordinary differential equations in the variables $T_s$, $G_{h,j}$, $T_{u,j}$, and $E$.

We assume the following values
for the system parameters: $K_{ps} = 0.2$, $T_h = 293.15$, $T_l = 288.15$, $T^0_s = 343.15$, $G_{hn} = 200$, $G_u = 150$, $Q_{max} = 0.7 N G_{hn}(T^0_s - T_h)$, $C_s = 2\cdot 10^6 N$, $t_h = 20$, $C_{u,j} = (1 + 0.348 j/N)\times 10^7$, $K_{pu} = 1$.
All the values of physical constants are in SI units.

The functions $f_j(t)$ are periodic with a period of one day, i.e., 86400 s. They start at $T_l$ at midnight, get increased from $T_l$ to $T_h$ at a pseudo-random time between 6 and 12 am, and switch back to $T_l$ at a pseudo-random time between 15 and 22 pm. The transitions are smooth, using the function $\operatorname{smoothStep}(\cdot, \cdot, \cdot)$ defined as
\begin{equation}
    \operatorname{smoothStep}(t, t_s, \Delta t) = \frac{1}{2}\left( \tanh \left( \frac{t - t_s}{\Delta t} \right) + 1 \right),
\end{equation}
with $\Delta t = 1$ s, while the function $\operatorname{sat}(\cdot, \cdot, \cdot)$ is a smooth saturation function defined as
\begin{eqnarray}
\operatorname{sat}(x, x_{min}, x_{max}) &=&  \frac{x_{max} + x_{min}}{2} \\
&+&  \frac{x_{max} - x_{min}}{2} \tanh \left( 2\frac{x - x_{min}}{x_{max} - x_{min}} - 1 \right). \nonumber
\end{eqnarray}
The initial conditions of the system at $t = 0$ (corresponding to 00:00, midnight) are $T_s = T^0_s$, $G_{h,j} = 0$, $T_{u,j} = 288.15$, $E = 0$. The integration interval is two days, i.e $t \in \left[ 0, 172800 \right]$.

We have used the 
fourth order  ESDIRK method denoted as
ESDIRK4(3)6L[2]SA  in  \cite{kennedy:2016} to compute
a reference single-rate solution with relative and absolute tolerance $10^{-9}, $ a single-rate solution with relative and absolute tolerance $10^{-5}$ and multi-rate solution with relative and absolute tolerance $10^{-5}$.
For the multi-rate algorithms, the parameter values   $\alpha_{max} = 1.2$, $\alpha_{min} = 0.5$, $\alpha = 0.9$, $\phi = 0.05$, and $\beta = 1$ were used; the associated  
dense output reconstruction was used for the interpolation operator $Q(\tau)$, both in the multi-rate procedure to interpolate the slow variables when integrating the fast variables, and when re-sampling the solution over a regular time grid with 10 s intervals.

\begin{figure}[tbph]
    \centering
    \includegraphics{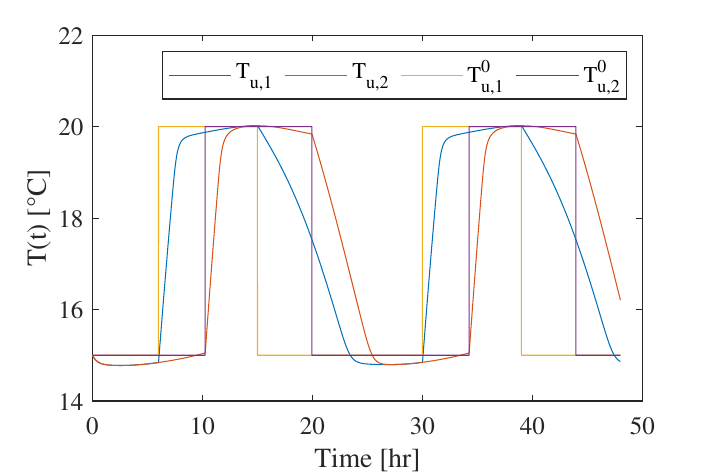}
    \caption{Temperatures and set points for units 1 and 2.}
    \label{fig:heating-unit-temp}
\end{figure}
\begin{figure}[tbph]
    \centering
    \includegraphics{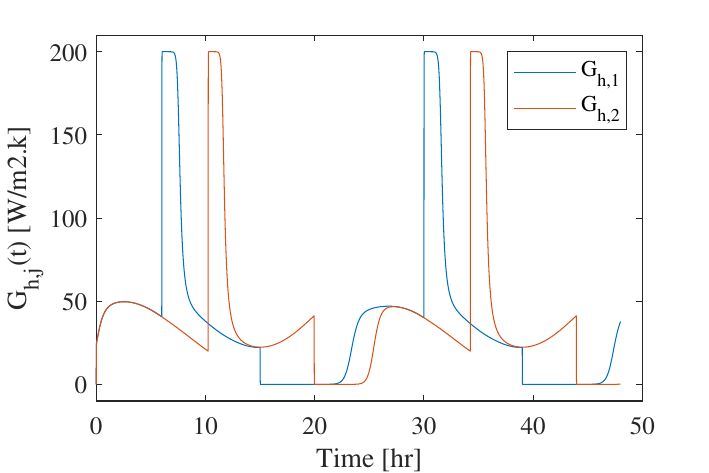}
    \caption{Conductance of fan-coils of units 1 and 2.}
    \label{fig:heating-actuators}
\end{figure}
\begin{figure}[tbph]
    \centering
    \includegraphics{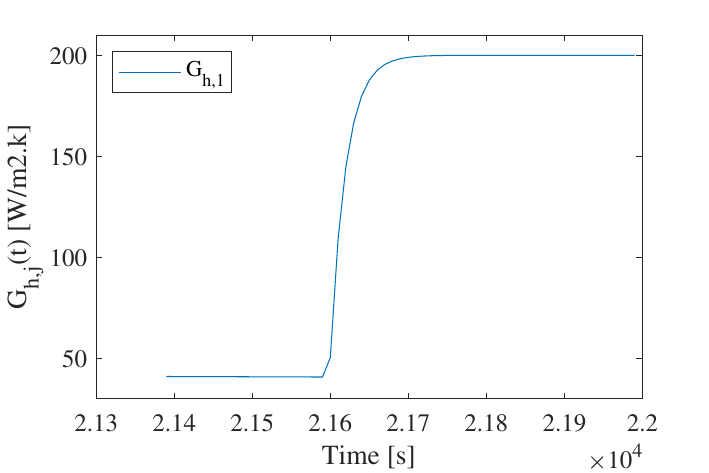}
    \caption{Detail of the first rising front of fan-coil 1.}
    \label{fig:heating-actuators-zoom}
\end{figure}

Figure \ref{fig:heating-unit-temp} shows the temperatures and the respective set points for units 1 and 2. The set point for the first unit $T^0_{u,1}$ is raised around 6:00 and reduced around 14:00, while the set point for the second unit $T^0_{u,1}$ is raised around 10:00 and reduced around 20:00. The unit temperatures $T_{u,1}$ and $T_{u,2}$ follow with some delay. Falling temperature transients are slower because they are only driven by the heat losses to the ambient.

Figure \ref{fig:heating-actuators} shows the fan-coil conductances $G_{h,1}$ and $G_{h,2}$. Initially, they are quite low, to keep the night temperature around $T_l.$ Subsequently, they very quickly reach the maximum value when the set point is raised, remain at the maximum for a while and then decrease once the unit temperatures approach the new set point value $T_h$, and drop sharply to zero when the set point is reduced to $T_l$ in the afternoon, increasing again during the night once the temperature has fallen below the set point. The dynamics of $G_{h,1}$ is much faster than the dynamics of the temperatures, see the detail of the first rising front of $G_{h,1}$ shown in Figure \ref{fig:heating-actuators-zoom}.

\begin{figure}[tbph]
    \centering
    \includegraphics{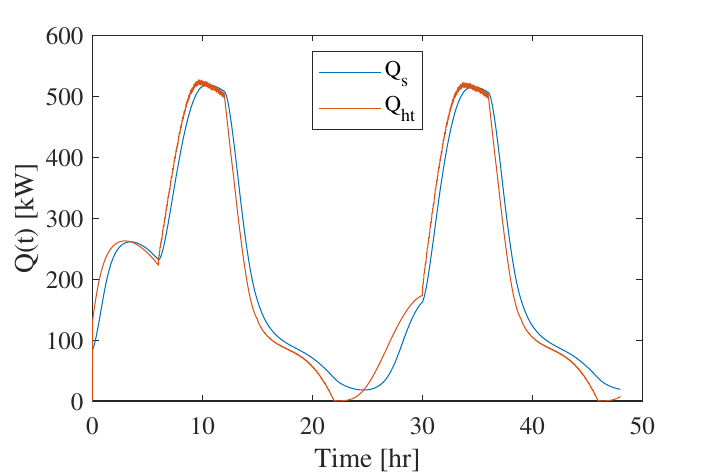}
    \caption{Input and output thermal power of the supply system.}
    \label{fig:heating-powers}
\end{figure}
\begin{figure}[tbph]
    \centering
    \includegraphics{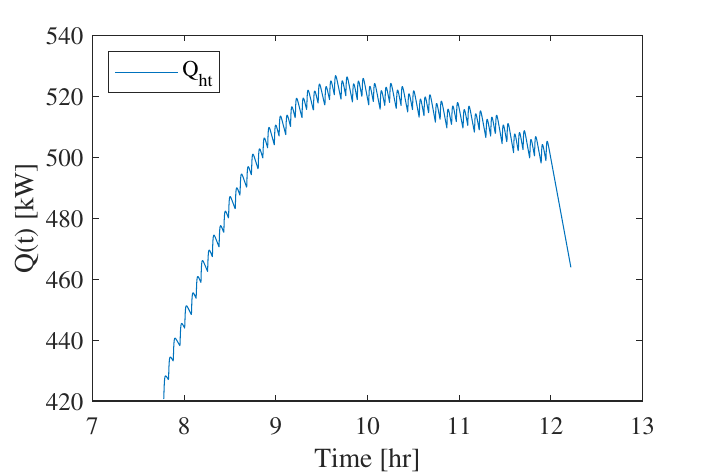}
    \caption{Zoom-in of the total output thermal power of the supply system between 7:45 and 12:15 on day 1.}
    \label{fig:heating-powers-zoom}
\end{figure}
\begin{figure}[tbph]
    \centering
    \includegraphics{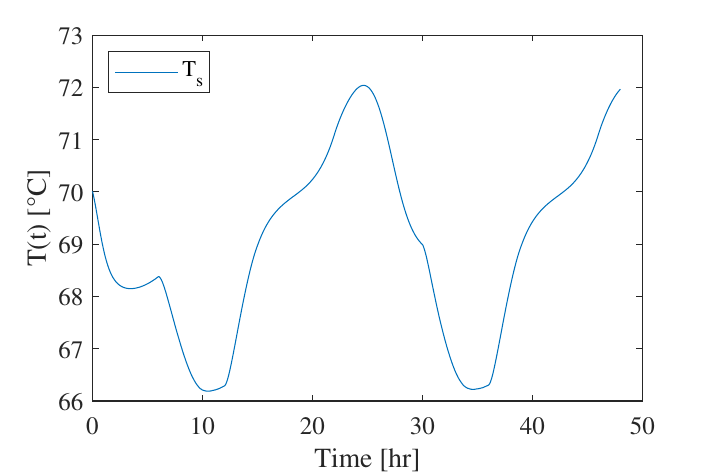}
    \caption{Heating supply temperature.}
    \label{fig:heating-supply-temp}
\end{figure}

Figure \ref{fig:heating-powers} shows the thermal power input of the supply system $Q_s$ and the total thermal power output to the unit fan coils $Q_{ht}$, with a zoom-in of $Q_{ht}$ between 7:45 and 12:15 shown in Figure \ref{fig:heating-powers-zoom}. The supply temperature $T_s$ is shown in Figure \ref{fig:heating-supply-temp}.

The reason why this system can benefit from multi-rate integration is twofold: the various occupants change their unit set points asynchronously, and the changes applied to one unit are very weakly coupled to the other units through the large thermal inertia of the supply system. Every time an occupant gets to its unit and raises the set point in the morning, or reduces it in the evening, a local fast transient is triggered on the unit's temperature control system output $u_j$, which requires relatively short time steps to simulate the transient of $G_{h,j}$. However, this action does not influence the other units directly, because they are insulated from each other, but rather only through the increased or decreased power consumption $Q_{ht}$, shown in Figure \ref{fig:heating-powers}. Even though $Q_{ht}$ shows fast changes, see Fig. \ref{fig:heating-powers-zoom}, corresponding to the individual unit heating systems being turned on or off, the large inertia of the supply system causes the supply temperature $T_s$ to remain quite smooth, as shown in Figure \ref{fig:heating-supply-temp}, so that the influence of these events on the heat exchanged by the other units becomes relevant only over much larger time intervals, even more so as $N$ grows.

Hence, every set point switching transient can be handled by a set of fast variables that only includes the local heater conductance $G_{h,j}$, the unit temperature $T_{u,j}$, the supply temperature $T_s$, and the consumed energy $E$, using interpolated values for the other unit variables, $G_{h,k}, T_{u,k}, k \neq j$, which only start getting influenced by the consequences of switching in other units after hundreds of seconds. This brings down the number of DOF required to manage these transients from $2N+2$ to $O(1)$.

\begin{figure}[tbph]
    \centering
    \includegraphics[width=\textwidth]{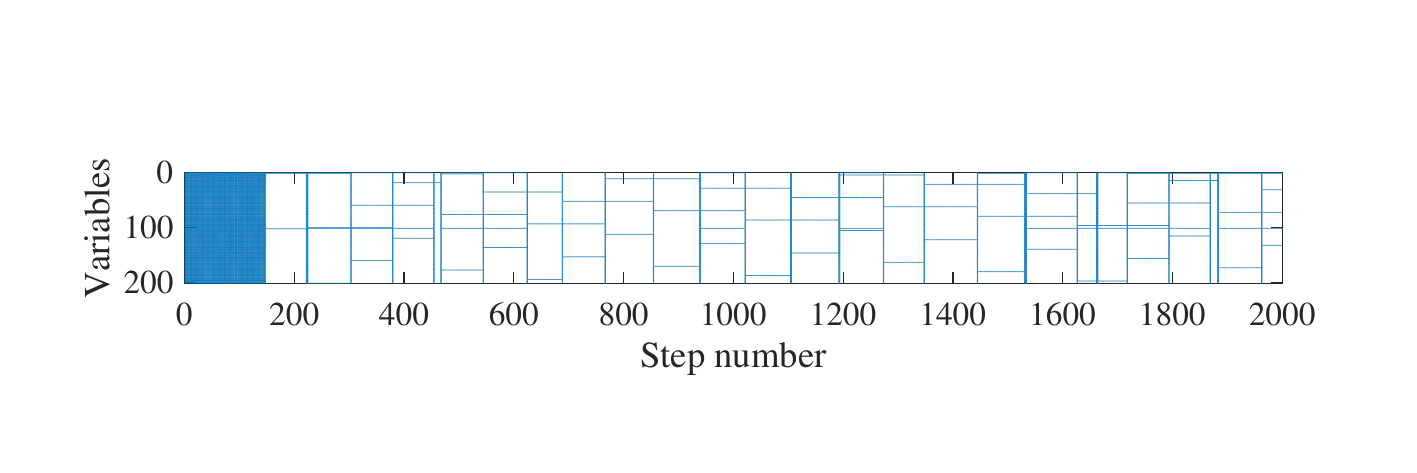}
    \includegraphics[width=\textwidth]{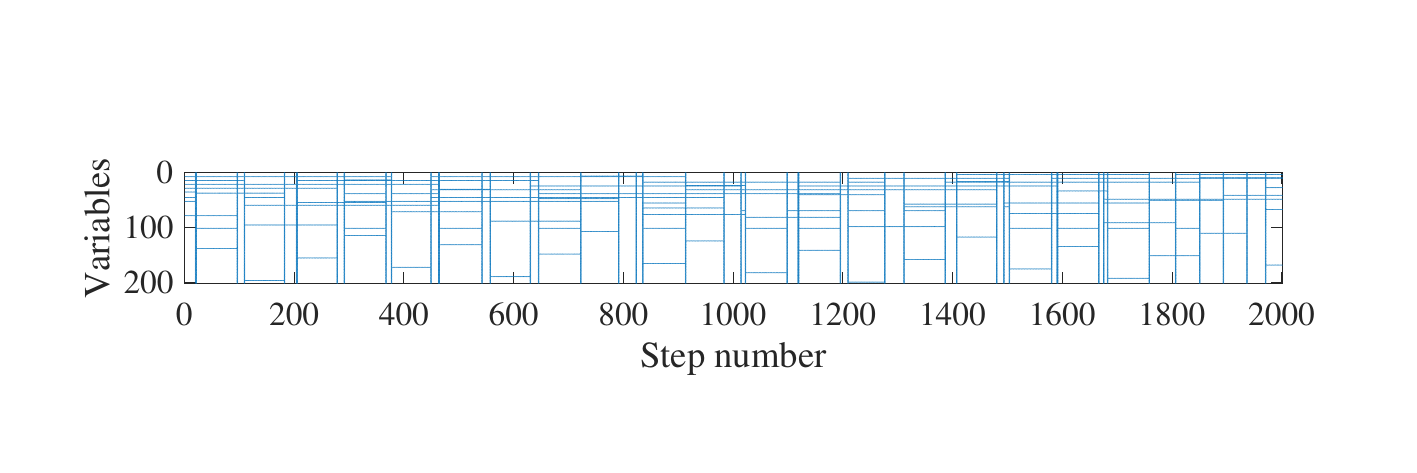}
    \includegraphics[width=\textwidth]{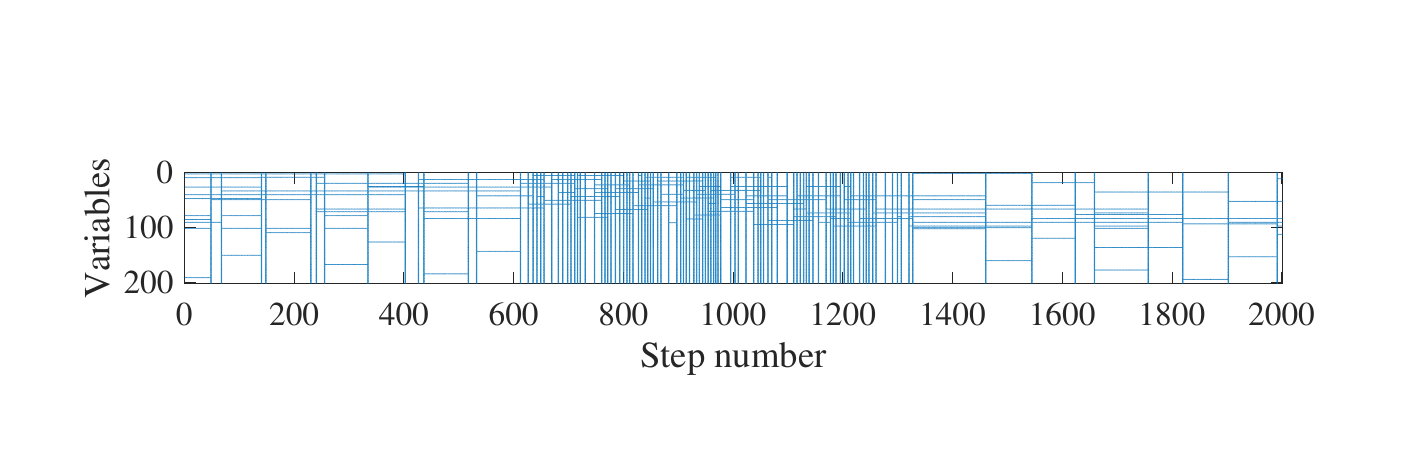}
    \caption{Activity diagrams of the multi-rate solution.}
    \label{fig:heating-activity}
\end{figure}
{   Figure \ref{fig:heating-activity}} shows three portions of the activity diagram, with rows corresponding to variables and columns to time steps. The top diagram corresponds to the beginning of the simulation, from 00:00 to 07:30. The first 145 steps to the left of the diagram correspond to the night transient from 00:00 to 06:00. During this period, no set points are changed, so there are no fast local transients taking place in the system, only the response of all the system variables to the slow sinusoidal variation of the external ambient temperature, which is handled by global time steps up to 1200 s long.

After 06:00, the unit temperature set points start being increased from $T_l$ to $T_h$, in a pseudo-random fashion, triggering many local fast transient. As a consequence, the average step length is dramatically reduced, and the remaining portion of the topmost activity diagram, spanning until 07:30, shows a very sparse pattern, where the two bordering lines at the top and bottom correspond to $T_u$ and $E$, which are always active, while the other shorter horizontal segments correspond to the transients of $G_{h,j}$ and $T_{u,j}$, which are handled by fast refinement steps; only a few global steps are necessary, shown in the figure as vertical blue bands).

This pattern continues until 12:00, e.g., see the the middle diagram of Figure \ref{fig:heating-activity}, covering the time interval from 8:45 to 10:00.

The bottom diagram of Figure \ref{fig:heating-activity}, which covers the interval from 11:30 to 15:30, shows the sparse patterns ending at 12:00, when set points stop being changed, followed by a pattern composed of mostly global steps until between 12:00 and 15:00, when the room temperatures mostly react to the slowly changing external temperature, after which set points start changing again, bringing the sparse pattern in again, and so on.

The most important result of this simulation is the final value of the energy consumption $E$, reported in {   Table \ref{tab:energy}} for the three solutions. The multi-rate result is affected by a larger error than the single-rate one, but it still has 5 correct significant digits, which is in agreement with the set relative and absolute tolerances $10^{-5}$, despite them being set on the local integration error and not on the cumulated error that affects the total energy consumption. 

\begin{table}[tbh]
    \centering
    \begin{tabular}{|c|c|}
    \hline    
                    & Total Consumption [MWh] \\
    \hline
       Reference   & 10.16868431299 \\
    \hline
       Single-Rate & 10.16868431246 \\
    \hline
       Multi-Rate  & 10.16822913061 \\
    \hline
    \end{tabular}
    \caption{Total energy consumption computed by the ESDIRK4(3)6L[2]SA method.}
    \label{tab:energy}
\end{table}

\begin{table}[tbh]
    \centering
    {  
    \begin{tabular}{|l|c|c|c|}
    \hline
                                                   & Single-Rate & Multi-Rate &  CVODE \\
       \hline
       Accepted global  steps                      & 27642       &  1106   & 42257 \\
       \hline
       Rejected global  steps (error test failure) & 5918        &  74     & 2470  \\
       \hline
       Accepted fast  steps                        &    //       &  31074  &   // \\
       \hline
       Rejected  fast steps                        &    //       &  5680   &   // \\
       \hline
       Global RHS calls                            &  1022880    &  73169  &  55042 \\
       \hline
       Global Jacobian computations                &  763        &  193    &   788  \\
       \hline
       Local RHS calls                            &    //        & 1056965 &   // \\
       \hline
       Local Jacobian computations                &    //        &  2546   &   // \\
       \hline
       Simulation time in s                       &    14.7       &  2.62  &  5.57 \\
       \hline
    \end{tabular}
    }
    \caption{Performance of single- and multi-rate ESDIRK4(3)6L[2]SA method.}
    \label{tab:heating}
\end{table}
The performance comparison between the single-rate and the multi-rate ESDIRK4(3)6L[2]SA method is shown in {   Table \ref{tab:heating}. As in the previous two cases, the number of total time steps required by the multi-rate method is only slightly larger than the number of steps of the single-rate method; in this case, the number of global steps is reduced by a factor 25 compared to the single-rate case. The speed-up factor of the multi-rate method over the single rate one is around 5.6.

In this case, unfortunately, the IDA solver did not produce satisfactory results in OpenModelica, due to some issues in the computation of the Jacobians. Therefore, the stiff CVODE solver was used instead, which is also implemented in SUNDIALS and is based on a BDF method. As shown in the last column of Table \ref{tab:heating}, CVODE is still faster than the single-rate implementation of ESDIRK4 in the GBODE solver, mainly due to the reduced number of right-hand-side computations. On the other hand, it turns out to be slower than the multi-rate ESDIKR4 method.

This already good result  can  be improved in the future by further optimizations of our code, both regarding the efficient solution of the global time step implicit equations and the efficient selective computation of the required parts of the right hand side for fast refinement states.}

 \section{Conclusions and future work}
 \label{sec:conclu} \indent
 {   We have presented a novel approach to  self-adjusting multi-rate methods, that allows to obtain an effective multi-rate implementation of a generic Runge-Kutta  (RK) method by
combining
a self-adjusting multi-rate technique with standard time step adaptation methods.}
Only a small percentage of the variables, associated with the
largest values of an error estimator are marked as fast variables. When the global step is sufficient to guarantee a given error tolerance for the slow variables, but not for the fast ones, the multi-rate procedure is employed to achieve uniform accuracy at reduced computational cost.  

We have  also derived a general  linear stability analysis, valid for  explicit RK,   DIRK and  ESDIRK multi-rate methods
with an arbitrary number of sub-steps for the active components. 
{   We have presented a thorough review of the literature on multi-rate stability analyses, which highlights the advantages of our approach. Furthermore, we have introduced a novel, physically motivated model problem that allows to assess the stability of multi-rate approaches in a context more relevant for realistic applications.}
 The   stability analysis has been performed on some examples of ERK and ESDIRK methods, highlighting the potential gains in efficiency and accuracy that can be achieved by the multi-rate approach  in the regime of weak coupling between fast and slow variables.
 
{  
The proposed multi-rate methods  have been  implemented efficiently in the framework of the OpenModelica software \cite{fritzson:2020}.
}
Three  numerical benchmarks have been considered, representing two idealized engineering systems
and a basic discretization of a nonlinear PDE typical of Computational Fluid  Dynamics applications.
{  
The results obtained demonstrate the   significant cost reductions with respect to the corresponding single-rate methods implemented the same framework and
a computational efficiency that starts to be comparable to that of much more mature single-rate solvers.}

In future work, we will further pursue the {   optimization of the present OpenModelica implementation.
In particular, the need for different  optimization strategies  in the multi-rate and single-rate case, highlighted by our numerical experiments, will be further investigated, in order to take full advantage of the multi-rate DOF reduction.
This will allow 
the effective application of multi-rate methods to larger scale PDE problems
as well as to more realistic engineering problems.  } 

\section*{Acknowledgements}
 {   We thank the reviewer and the Associate Editor for their constructive remarks, which greatly helped in improving the first version of this paper.} B.B. carried out part of this research during a sabbatical at Politecnico di Milano in 2022, which was partially  supported by the Visiting Faculty Program of Politecnico di Milano. S.F.G, M.G.M and L.B. have also received support through the PID2021-123153OB-C21 project of the Science and Innovation Ministry of the Spanish government.

\appendix
 \section{Butcher tableaux of  RK methods}
 \label{sec:butcher} \indent
For completeness, we report in this Appendix  the Butcher tableaux of the methods considered in the 
stability analysis, along with the
coefficients of their corresponding dense output interpolators.
  Table \ref{rk4tab} contains the coefficients defining the classical fourth order explicit RK method.
Table \ref{rkzennaro} contains the coefficients defining the fourth order RK method with optimal continuous output
 introduced in \cite{owren:1992}. 
The matrix $B^*$ containing the coefficients for the dense output formula \eqref{eq:dense1} associated to this method are given in Table \ref{rkzennaro_bstar}. Setting then
\begin{eqnarray}
  \gamma&=&0.43586652150845899941601945, \nonumber \\  
  c_{3}&=&\frac35,  \ \ \ \ \bar a=1-6\gamma+6\gamma^2 \nonumber \\
  a_{32}&=&c_3\frac{c_3-2\gamma}{4\gamma} \ \ \  a_{31}=c_3-a_{32}-\gamma  \\
  b_2&=&\frac{-2+3c_3+6\gamma(1-c_3)}{12\gamma(c_3-2\gamma)} \ \ \ 
  b_3=\frac{\bar a}{3c_3(c_3-2\gamma)} \nonumber \\
  b_1&=&1-b_2-b_3-\gamma, \nonumber
\end{eqnarray}
we report in Table \ref{esdirk3} the coefficients of the third order ESDIRK method denoted as
ESDIRK3(2)4L[2]SA  in  \cite{kennedy:2016}.
The corresponding matrix $B^*$ containing the coefficients for the dense output formula \eqref{eq:dense1} are given in Table \ref{esdirk3_bstar}.

\small
\begin{table}
 \begin{center}  
\begin{tabular}{c|cccc} 
0 & 0 & 0& 0&0\\
 &  & & &\\
$ \frac 12 $ & $ \frac 12$& 0& 0& 0\\
 &  & & &\\
$\frac 12$ & $0$&$\frac 12$ & 0& 0\\
 &  & & &\\
 $1$ & $0$& $0$& 1 &0 \\
  &  & & &\\
\hline 
 &  & & &\\
 &  $\frac 16$ & $\frac 26$ & $ \frac 26$ & $\frac 16 $\ \
\end{tabular}
\caption{Butcher tableaux of the classical fourth order ERK method. \label{rk4tab} }
 \end{center}
 \end{table} 
\normalsize

\small
\begin{table}
 \begin{center}  
\begin{tabular}{c|cccccc} 
 0                &   0             &0  &0   &0 &0 &0\\
                   &                   &    &    &   &   & \\
$ \frac 16 $ & $ \frac 16$& 0 & 0 & 0 &0 &0 \\
 &  & & & & & \\
$\frac {11}{37}$ & $\frac {44}{1369}$& $\frac {363}{1369}$ &  0 & 0& 0 & 0\\
 &  & & & & &\\
 $\frac {11}{17}$  & $\frac{3388}{4913}$& $-\frac{8349}{4913}$& $-\frac{8140}{4913}$ &  0 & 0 & 0\\
  &  & & & & &\\
$\frac {13}{15}$  & $-\frac{36764}{408375}$ & $\frac{767}{1125}$& $-\frac{32708}{136125}$  &$\frac{210392}{408375}$ &0 &0  \\
  &  & & & & &\\
$1$  & $\frac {1697}{18876}$ & 0 &  $\frac{50653}{116160}$ & $\frac{299693}{1626240}$ &
$\frac{3375}{11648} $ & 0\\
  &  & & & & &\\
\hline 
 &  & & & & &\\
 &  $\frac {101}{363}$ & $0$ & $ -\frac {1369}{14520}$ & $\frac {11849}{14520}$ &0 &0 \\
\end{tabular}
\caption{Butcher tableaux of optimal explicit RK method with continuous output. \label{rkzennaro} }
 \end{center}
 \end{table} 
 \normalsize

\begin{table}
 \begin{center}  
\begin{tabular}{|c|c|c|c|c|} 
\hline
 $b^*_{ij}$     &   $j=1$    &$j=2$  & $j=3$ &$j=4$ \\
  \hline
  & & & & \\
 $i=1$ & $  1 $ & $ -\frac{104217}{37466}$ & 
 $\frac{1806901}{618189}$  & $ -\frac{866577}{824252} $   \\
  & & & & \\
\hline
$i=2$  &$ 0$ & $0$& $0$ &  0  \\
\hline
 & & & & \\
 $i=3$ & $0$  & $\frac{861101}{230560}$ & $-\frac{2178079}{380424}$& $\frac{12308679}{5072320 }$  \\
  & & & & \\
 \hline  
  & & & & \\
$i=4$ &$0$  & $-\frac{63869}{293440}$ & $\frac{6244423}{5325936}$& $-\frac{7816583}{10144640}$  \\
 & & & & \\
\hline 
 & & & & \\
$i=5$ & $0$  & $-\frac {1522125}{762944}$ & $\frac {982125}{190736}$  &  $-\frac{624375}{217984}$   \\
 & & & & \\
\hline  
 & & & & \\
$i=6$ & $0$&  $\frac{165}{131}$  & $-\frac {461}{131}$ &  $\frac {296}{131}$   \\
 & & & & \\
\hline  
\end{tabular}
\caption{Matrix of continuous output coefficients $b^*_{ij}$ for explicit RK method of 
\cite{owren:1992}.  \label{rkzennaro_bstar}}
 \end{center}
 \end{table}

\begin{table}

 \begin{center}  
\begin{tabular}{c|cccc} 
0 & 0 & 0& 0&0\\
 &  & & &\\
$ 2\gamma $ & $ \gamma $& $ \gamma $ & 0& 0\\
 &  & & &\\
$\frac 35$ & $a_{31}$ 
& $a_{32}$ & $ \gamma $ & 0\\
 &  & & &\\
 $ 1$ & $b_1$& $b_2$& $b_3$ &$ \gamma $ \\
  &  & & &\\
\hline 
 &  & & &\\
 &  $b_1$ & $b_2$ & $b_3$ & $\gamma $\ \
\end{tabular}
\caption{Butcher tableaux of the ESDIRK3(2)4L[2]SA  method.\label{esdirk3} }
 \end{center}
 \end{table}

\begin{table}
 \begin{center}    
\begin{tabular}{|c|c|c|c|} 
\hline
 $b^*_{ij}$     &   $j=1$    &$j=2$  & $j=3$  \\
  \hline
  & & & \\
 $i=1$ & $ \frac{6071615849858}{5506968783323}  $ & 
 $ -\frac{9135504192562}{5563158936341}$ & 
 $\frac{5884850621193}{8091909798020}$    \\
  & & &  \\
\hline
 & & &  \\
$i=2$ & $ \frac{24823866123060}{14064067831369}  $ & 
 $ -\frac{184358657789355}{34679930461469}$ & 
 $\frac{40093531604824}{13565043189019}$    \\
  & & &  \\
\hline
 & & &  \\
 $i=3$ & $ -\frac{4639021340861}{5641321412596}  $ & 
 $ -\frac{36951656213070}{8103384546449}$ & 
 $-\frac{9445293799577}{3414897167914}$    \\
  & & &  \\
 \hline  
  & & & \\
  $i=4$ & $ -\frac{4782987747279}{4575882152666}  $ & 
 $ \frac{22547150295437}{9402010570133}$ & 
 $-\frac{8621837051676}{9402290144509}$    \\
 & & &  \\
\hline 
\end{tabular}
\caption{Matrix of continuous output coefficients $b^*_{ij}$ for the ESDIRK3(2)4L[2]SA  method. \label{esdirk3_bstar} }
 \end{center}
 \end{table}

\begin{table}
 \begin{center}  
\begin{tabular}{c|cccccc} 
 0                &   0             &0  &0   &0 &0 &0\\
                   &                   &    &    &   &   & \\
$2\gamma $ & $\gamma $& $\gamma $ & 0 & 0 &0 &0 \\
 &  & & & & & \\
$c_3$ & $a_{31}$& $a_{32}$ &  $\gamma $ & 0& 0 & 0\\
 &  & & & & &\\
 $c_4$  & $a_{41}$& $a_{42}$& $a_{43}$ &  $\gamma $ & 0 & 0\\
  &  & & & & &\\
$c_5$  & $a_{51}$ & $a_{52}$& $a_{53}$  & $a_{54}$  &$\gamma $ &0  \\
  &  & & & & &\\
$1$ &  $b_1$ & $b_2$  & $b_3$  & $b_4$ & $b_5 $ &$ \gamma$\\
  &  & & & & &\\
\hline 
 &  & & & & &\\
 &  $b_1$ & $b_2$  & $b_3$  & $b_4$ & $b_5 $ &$ \gamma$ \\
\end{tabular}
\caption{Butcher tableaux of the ESDIRK4(3)6L[2]SA method. \label{esdirk4} }
 \end{center}
 \end{table} 
 
 \pagebreak

  Table  \ref{esdirk4} contains instead the coefficients defining the 
   fourth order  ESDIRK method denoted as
ESDIRK4(3)6L[2]SA  in  \cite{kennedy:2016}, 
  where we have set instead 

\begin{eqnarray}
\gamma&=&1/4 c_3=\frac{2-\sqrt{2}}{4} \ \ \ \ 
c_4=\frac 58 \ \  c_5=\frac{26}{25} \nonumber \\
  a_{32}&=&\frac{1-\sqrt{2}}{8} \ \ 
 a_{31}=c_3-a_{32}-\gamma \nonumber \\
 a_{42}&=&\frac{5-7\sqrt{2}}{64} \ \ \ \   a_{43}=7\frac{1+\sqrt{2}}{32} \ \ \ \
 a_{41}=c_4-a_{42}-a_{43}-\gamma \nonumber \\
 a_{52}&=&\frac{-13796-54539\sqrt{2}}{125000} \ \ 
 a_{53}=\frac{506605+132109\sqrt{2}}{437500} \\
a_{54}&=&166\frac{-97+376\sqrt{2}}{109375} \ \ 
 a_{51}=c_5-a_{52}-a_{53}-a_{54}-\gamma \nonumber \\
   b_2&=&\frac{1181-987\sqrt{2}}{13782} \ \  \ \ 
 b_3=47\frac{-267+1783\sqrt{2}}{273343} \nonumber \\
 b_4&=&-16\frac{-22922+3525\sqrt{2}}{571953} \ \ \ \
 b_5=-15625\frac{97+376\sqrt{2}}{90749876} 
 \nonumber \\
 b_1&=&1-b_2-b_3-b_4-b_5-\gamma.\nonumber 
 \end{eqnarray}

The corresponding matrix $B^*$ containing the coefficients for the dense output formula \eqref{eq:dense1} are given in Table \ref{esdirk4_bstar}.

\pagebreak

\begin{table}
 \begin{center}    
\begin{tabular}{|c|c|c|c|c|} 
\hline
 $b^*_{ij}$     &   $j=1$    &$j=2$  & $j=3$ &$j=4$ \\
  \hline
  & & & & \\
 $i=1$ & $  \frac{11963910384665}{12483345430363} $ & $ 
 -\frac{69996760330788}{18526599551455}$ & 
 $\frac{32473635429419}{7030701510665}$  & $ -\frac{14668528638623}{8083464301755} $   \\
  & & & & \\
\hline
& & & & \\
 $i=2$ & $  \frac{11963910384665}{12483345430363} $ & $ 
 -\frac{69996760330788}{18526599551455}$ & 
 $\frac{32473635429419}{7030701510665}$  & $ -\frac{14668528638623}{8083464301755} $   \\
& & & & \\
\hline
 & & & & \\
 $i=3$ & $ -\frac{28603264624}{1970169629981}$  & $\frac{102610171905103}{26266659717953}$ & $-\frac{38866317253841}{6249835826165}$& 
 $\frac{21103455885091}{7774428730952}$  \\
  & & & & \\
 \hline  
  & & & & \\
$i=4$ & $-\frac{3524425447183}{2683177070205}$ & $\frac{74957623907620}{12279805097313}$& $-\frac{26705717223886}{4265677133337 }$
&$\frac{30155591475533}{15293695940061}$\\
 & & & & \\
\hline 
 & & & & \\
$i=5$ & $ -\frac{17173522440186}{10195024317061}$  & 
$\frac{113853199235633}{9983266320290}$ & $-\frac{121105382143155}{6658412667527}$  &  $\frac{119853375102088}{14336240079991}$   \\
 & & & & \\
\hline  
 & & & & \\
$i=6$ & $ \frac{27308879169709}{13030500014233 }$&  $-\frac{84229392543950}{6077740599399}$  & $\frac {1102028547503824}{51424476870755}$ &  $-\frac {63602213973224}{6753880425717}$   \\
 & & & & \\
\hline  
\end{tabular}

\caption{Matrix of continous output coefficients $b^*_{ij}$ for the ESDIRK4(3)6L[2]SA method. \label{esdirk4_bstar} }
 \end{center}
 \end{table}

\bibliographystyle{plain}
\bibliography{multirate_2025}

\end{document}